\def\versiondate{30 April 1999}
\input math.macros
\input Ref.macros

%\proofmodetrue
%\proofmodefalse
%\leftsectionheadtrue
\checkdefinedreferencetrue
%\continuousnumberingtrue
\continuousfigurenumberingtrue
\theoremcountingtrue
\sectionnumberstrue
%\figuresectionnumberstrue
\forwardreferencetrue
%\lefteqnumberstrue
\hyperstrue
\initialeqmacro

\def\st{:\,}
\def\ev#1{{\cal #1}} % Events
\def\A{\ev A}
\def\D{\ev D}
\def\verts{{\ss V}}
\def\edges{{\ss E}}
\def\xor{\bigtriangleup}
\def\wtA{{\ev A}'_o}
\def\mnote#1{\ifproofmode{\bf [#1] }\fi}
\def\H{{\Bbb H}}
\def\ins{\Pi}   %% Insertion
\def\Aut{{\rm Aut}}
\def\assumptions{Let $G$ be a graph with a transitive unimodular
                 closed automorphism group $\Gamma\subset\Aut(G)$}
\def\assumptionsandinf{Let $G$ be an infinite graph with a transitive unimodular
                 closed automorphism group $\Gamma\subset\Aut(G)$}
\def\genassumptions{Let $\Gamma$ be a transitive closed subgroup of the
                 automorphism group of a graph $G$}
\def\SS{{\cal S}}   %% shift 
\def\setminu{-}
\def\Ph{\widehat \P}   %% probability measure on very big space
% \def\Omh{\widehat \Omega}   %% big space
   %% invariant infinite measure on very big space
\def\freq{\alpha}
\def\freqfn{{\rm freq}}
\def\dist{{\rm dist}}
\def\bd{\partial}
\def\U{{\cal U}}
\def\H{{\cal H}}
\def\F{{\cal F}}
\def\eps{\epsilon}

\def\id{{\rm id}}

   %% shift and multiply 
\def\wreath{\wr}
\font\frak=eufm10   %% or  eufb10
\def\fo{{\hbox{\frak F}}}

\def\law{{\bf Q}_p}  %% law of $\eta$

\def\BLPSgip{Benjamini, Lyons, Peres, and Schramm (1999)%
\def\BLPSgip{BLPS (1999)}}
\def\BLS{Benjamini, Lyons, and Schramm (1998)}%

\def\firstheader{\eightpoint\ss\underbar{\raise2pt\line 
    {To appear in {\it Ann.\ Probab.}\hfil Version of \versiondate}}}

\vglue20pt

\beginniceheadline

%\ifproofmode \relax \else\head{} {Version of \versiondate}\fi 
%\vglue20pt

\title{Indistinguishability of Percolation Clusters}

\author{Russell Lyons and Oded Schramm}

\abstract{
We show that when percolation produces infinitely many infinite clusters on
a Cayley graph, one cannot distinguish the clusters from each
other by any invariantly defined property. This implies that
uniqueness of the infinite cluster is equivalent to non-decay of connectivity
(a.k.a.\ long-range order).
We then derive applications concerning uniqueness in
Kazhdan groups and in wreath products,
and inequalities for $p_u$.
}

\bottom{Primary 82B43, % Percolation
60B99. % probability thy on algebraic and topological structures
 Secondary
60K35, % Interacting random processes; statistical mechanics type models; percolation theory
60D05. % geometric probability, random sets
} 
{Finite energy, Cayley graph, group, Kazhdan, wreath product, uniqueness,
connectivity, transitive, nonamenable. 
}
{Research partially supported by a Varon Visiting Professorship at
the Weizmann Institute of Science and NSF grant DMS-9802663 (Lyons) and
the Sam and Ayala Zacks Professorial Chair (Schramm).}

\bsection {Introduction}{s.intro}

Grimmett and Newman (1990) showed that if $T$ is a
regular tree of sufficiently high degree, then there are
$p\in(0,1)$ such that Bernoulli($p$) percolation on
$T\times \Z$ has infinitely many infinite components a.s.
Benjamini and Schramm (1996) conjectured that the same
is true for any Cayley graph of any finitely generated nonamenable group.
(A finitely generated group $\Gamma$ is nonamenable
iff its Cayley graph satisfies 
$\inf_K |\partial K|/|K|>0$, where $K$ runs over the finite
nonempty vertex sets. See \ref s.back/ for all other definitions.)
This conjecture has been verified for
planar Cayley graphs of high genus by Lalley (1998) and for all
planar lattices in the hyperbolic plane
by Benjamini and Schramm (in preparation).

The present paper is concerned with the percolation phase where there are
multiple infinite clusters. We show that under quite
general assumptions, the infinite clusters are indistinguishable from each
other by any invariantly defined property.

Let $G$ be a Cayley graph of a finitely generated group
(more generally, $G$ can be a transitive unimodular graph).

\procl t.main
Consider Bernoulli
bond percolation on $G$ with some survival parameter $p\in(0,1)$,
and let $\ev A$ be a Borel measurable set of subgraphs of $G$.
Assume that $\ev A$ is invariant under the automorphism group of $G$. 
Then either a.s.\ all infinite percolation components are in $\ev A$,
or a.s.\ they are all outside of $\ev A$.
\endprocl

For example, $\ev A$ might be the collection of all transient
subgraphs of $G$, or the collection of all subgraphs that have
a given asymptotic rate of growth, or the collection of all subgraphs
that have no vertex of degree $5$.

If $\ev A$ is the collection of all transient subgraphs
of $G$, then this shows that almost surely, either all infinite
clusters of $\omega$ are
transient [meaning that simple random walk on them is transient],
or all clusters are recurrent.
In fact, as shown in \ref p.trans/, if $G$ is nonamenable,
then a.s.\ all infinite clusters are transient if Bernoulli percolation
produces more than one infinite component.  In \BLS, it is shown
that the same is true if Bernoulli percolation produces a single
infinite component.

\ref t.main/ is a particular case of \ref t.cerg/ below,
which applies also to some non-Bernoulli percolation processes
and to more general $\ev A$.

A collection of subgraphs $\ev A$ in $G$ is {\bf increasing}
if whenever $H\in\ev A$ and $H\subset H'\subset G$, we have $H'\in\A$.
H\"aggstr\"om and Peres (1999) have proved that for
increasing $\A$, \ref t.main/ holds, except possibly for a single
value of $p\in(0,1)$. They have also proved:

\procl t.uniq \procname{Uniqueness Monotonicity}
Let $p_1 < p_2$ and $\P_i$ ($i=1,2$) be the corresponding Bernoulli($p_i$) bond
percolation processes on $G$.
If there is a unique infinite cluster $\P_1$-a.s., then
there is a unique infinite cluster $\P_2$-a.s.  Furthermore, in the
standard coupling of Bernoulli percolation processes, if there
exists an infinite cluster $\P_1$-a.s., then a.s.\ every
infinite $\P_2$-cluster contains an infinite $\P_1$-cluster.
\endprocl

Here, we refer to the standard coupling of Bernoulli($p$) percolation for
all $p$, where each edge $e \in \edges$ is assigned an independent uniform
$[0, 1]$ random variable $U(e)$ and the edges where $U(e) \le p$ are
retained for Bernoulli($p$) percolation. 

The first part of this theorem partially answers a question of Benjamini
and Schramm (1996); the full answer, removing the assumption of
unimodularity, has been provided by Schonmann (1999a). 
It follows from \ref t.uniq/ that the set of $p$ such that
Bernoulli($p$) bond percolation on $G$ has more than one infinite
cluster a.s.\ is an interval, called the {\bf nonuniqueness phase}.
It is well known that in the nonuniqueness phase, the number of infinite
clusters is a.s.\ infinite (see Newman and Schulman 1981). 
The interval of $p$ such that
Bernoulli($p$) percolation on $G$ produces a single infinite cluster
a.s.\ is called the {\bf uniqueness phase}.  The infimum of the
$p$ in the uniqueness phase is denoted $p_u=p_u(G)$.

To illustrate the usefulness of \ref t.main/, we now show that
it implies \ref t.uniq/.

\proofof t.uniq
Suppose that there exists an infinite cluster $\P_1$-a.s.
Let $\omega$ be the open subgraph of the $\P_2$ process and let $\eta$ be
an independent Bernoulli($p_1/p_2$) percolation process.
Thus, $\omega \cap \eta$ has the law of $\P_1$ and, in fact, $(\omega \cap
\eta, \omega)$ has the same law as the standard coupling of $\P_1$ and
$\P_2$.
By assumption, $\omega \cap \eta$ has an infinite cluster a.s.
Thus, for some cluster $C$ of $\omega$, we have $C \cap \eta$ is infinite with
positive probability, hence, by Kolmogorov's 0-1 law, with probability 1.
By \ref t.main/, this holds for every cluster $C$ of $\omega$.
\Qed 

\comment{
Because of the importance of \ref t.uniq/, we give
another short proof of \ref t.uniq/ (that does not use
\ref t.main/) in \ref r.easy/ below.
}

Our main application of \ref t.main/ is to prove in \ref s.decay/
a sufficient condition for having a unique infinite cluster.
Define $\tau(x, y)$ to be
the probability that $x$ and $y$ are in the same cluster.  
The process is said to exhibit {\bf connectivity decay}
if 
$$
\inf\bigl\{ \tau(x,y)\st x, y \in \verts\bigr\} = 0
\,.
$$
When a process does not have connectivity decay, 
it is said to exhibit {\bf long-range order}.
We show that for Bernoulli percolation (and some other percolation
processes), nonuniqueness of the infinite cluster is equivalent
to connectivity decay.
One direction is easy, namely, if with positive probability
there is a unique infinite percolation cluster,
then there is long-range order.  This is an immediate consequence
of Harris's inequality
(that increasing events are positively correlated).
Although the other direction seems intuitively obvious as well,
its proof is surprisingly difficult
(and it is still open for nonunimodular transitive graphs).

We note that connectivity decay in the nonuniqueness phase
easily implies Thm.~3.2 of Schonmann (1999a) for
the case of unimodular transitive graphs.
(Our result also deals with percolation at
$p_u$ and extends Schonmann's theorem to more general percolation
processes, but his result extends to Bernoulli percolation on nonunimodular
transitive graphs.)

\procl r.directions
Fix a base vertex $o\in\verts$ in $G$, and consider
Bernoulli percolation in the nonuniqueness phase.  
Although $\inf\bigl\{ \tau(o,v)\st{v\in\verts}\bigr\}=0$, it may happen that
$$
\lim_{d\to\infty} \sup\bigl\{\tau(o,v)\st\dist(o,v)>d\bigr\}>0
%\limsup_{\dist(o,v)\to\infty} \tau(o,v)>0
\,.
$$
For example, take the Cayley graph of the free product
$\Z^2*\Z_2$ with the usual generating set, i.e., the Cayley graph
corresponding to the presentation
$\langle a, b, c \mid ab=ba, c = c^{-1} \rangle$.
We have $p_u=1$,
because the removal of any edge of the Cayley graph 
corresponding to the generator $c$ disconnects the graph.
On the other hand, for $p > p_c(\Z^2)$, we have that all
pairs of vertices in the same $\Z^2$ fiber
(i.e., pairs $g, h$ such that $g^{-1} h$ is in the group generated
by $a$ and $b$) are
in the same cluster with probability bounded away from 0.
\endprocl

Benjamini and Schramm (1996) asked for which Cayley graphs $G$ one has
$p_u(G) < 1$. 
It is easy to show that $p_u(G) = 1$ when $G$ has more than one end.
Babson and Benjamini (1999) showed that $p_u(G) < 1$ when $G$ is a
finitely presented group with one end.
H\"aggstr\"om, Peres and Schonmann (1998) have shown that
$p_u(G)<1$ if $G$ is a Cartesian product of infinite transitive graphs.
Here, we show that this is also true of certain other classes of groups
with one end, namely infinite Kazhdan
groups (\ref c.kazhdan/) and wreath products (\ref c.wreath/).
Presumably, all quasi-transitive (infinite) graphs with one end have $p_u < 1$.

As another consequence of indistinguishability,
we prove in \ref t.puineq/ the inequalities
$p_u(G\times H')\geq p_u(G\times H)$ and, in particular,
$p_u(G) \geq p_u(G\times H)$, for Cayley graphs
(or, more generally, unimodular transitive graphs)
$G,H',H$ such that $G$ is infinite and $H'\subset H$.

Several other uses of cluster indistinguishability appear in \BLS.

Crucial techniques for our proofs are the Mass-Transport Principle and
stationarity of delayed simple random walk, both explained below. These
techniques were introduced in the study of percolation by H\"aggstr\"om (1997).
In \ref s.ergdrw/, we prove an ergodicity property for delayed random walk.

\bsection {Background}{s.back}

Let $\Gamma$ be a finitely generated group and $S$ a finite generating set for
$\Gamma$.  Then the (right) {\bf Cayley} graph $G=G(\Gamma,S)$ is the
graph with vertices $\verts(G):=\Gamma$ and edges
$\edges(G):=\bigl\{[v,vs]\st v\in\Gamma,\, s\in S\bigr\}$.

Now suppose that $G=\big(\verts(G),\edges(G)\big)$ is any graph,
not necessarily a Cayley graph.
An {\bf automorphism} of $G$ is a bijection of $\verts(G)$ with
itself that preserves adjacency; $\Aut(G)$ denotes the group of all
automorphisms of $G$ with the topology of pointwise convergence.
If $\Gamma\subset\Aut(G)$, we say that $\Gamma$ is (vertex) {\bf
transitive} if for every $u, v\in\verts(G)$, there is a $\gamma\in\Gamma$
with $\gamma u=v$.  
The graph $G$ is {\bf transitive} if $\Aut(G)$ is transitive.
The graph $G$ is {\bf quasi-transitive} if $\verts(G)/\Aut(G)$ is finite;
that is,
there is a finite set of vertices $V_0$ such that 
$\verts(G)=\bigl\{\gamma v\st \gamma\in\Aut(G),\, v\in V_0\bigr\}$.
Note that any finitely generated group acts transitively on any 
of its Cayley graphs by the automorphisms $\gamma : v \mapsto \gamma v$.

It may seem that the most natural class of graphs
on which percolation should be studied is the class of transitive
(or quasi-transitive) locally finite graphs. 
At first sight, one might suspect that any theorem about percolation
on Cayley graphs should hold for transitive graphs.
However, somewhat surprisingly, this impression is not correct.
It turns out that theorems about percolation on Cayley graphs
``always" generalize to unimodular transitive graphs (to be defined
shortly), but nonunimodular transitive graphs are quite different.

Recall that on every closed subgroup $\Gamma\subset \Aut(G)$, there is
a unique (up to a constant scaling factor) Borel measure that, for every
$\gamma\in\Gamma$, is invariant under left multiplication by $\gamma$; this
measure is called (left) {\bf Haar measure}.
The group $\Gamma$ is {\bf unimodular} if Haar measure is also invariant
under right multiplication.
For example, when $\Gamma$ is finitely generated and $G$ is the (right)
Cayley graph of $\Gamma$, on which $\Gamma$ acts by left multiplication,
then $\Gamma\subset\Aut(G)$
is (obviously) closed, unimodular, and transitive.  The Haar measure
in this case is (a constant times) counting measure.
A [quasi-]transitive graph $G$ is said to be {\bf unimodular} if
$\Aut(G)$ is unimodular.  By \BLPSgip{} (hereinafter referred to as
\BLPSgip), a transitive graph
is unimodular iff there is some unimodular transitive closed subgroup
of $\Aut(G)$.
It is not hard to show that a transitive closed subgroup
$\Gamma\subset\Aut(G)$ is unimodular iff for all $x, y \in \verts(G)$, we have
$$
\big|\big\{z \in \verts(G) \st \texists {\gamma \in \Gamma} \gamma x = x
\hbox{ and } \gamma y = z \big\}\big|
=
\big|\big\{z \in \verts(G) \st \texists {\gamma \in \Gamma} \gamma y = y
\hbox{ and } \gamma x = z \big\}\big|
$$
(see Trofimov (1985)).

Here, unimodularity will be used only in the following form:

\procl t.mtp \procname{Unimodular Mass-Transport Principle}
\assumptions.
Let $o\in\verts(G)$ be an
arbitrary base point.  Suppose that
$\phi:\verts(G)\times\verts(G)\to [0,\infty]$
is invariant under the diagonal action of $\Gamma$.
Then
$$
\sum_{v\in\verts(G)} \phi(o,v) = \sum_{v\in\verts(G)} \phi(v,o)
\,.
\label e.mtp
$$
\endprocl

See \BLPSgip\ for a discussion of this principle and for a proof.
In fact,  $\Gamma$ is unimodular iff \ref e.mtp/ holds for
every such $\phi$.  Also see \BLPSgip\ for further discussion
of the relevance of unimodularity to percolation.

We now discuss some graph-theoretic terminology.
When there is an edge in $G$ joining vertices $u,v$, we say that $u$ and
$v$ are {\bf adjacent} and write $u\sim v$.
We always assume that the number of vertices adjacent to any given vertex
is finite. 
The {\bf degree} $\deg v =\deg_G v $ of a vertex $v\in\verts(G)$ is the
number of edges incident with it.
A {\bf tree} is a connected graph with no cycles.
A {\bf forest} is a graph whose connected components are
trees.
The {\bf distance} between two vertices $v,u\in\verts(G)$
is denoted by $\dist(v,u)=\dist_G(v,u)$, and is the least number of edges
of a path in $G$ connecting $v$ and $u$. 
Given a set of vertices $K\subset \verts(G)$, we let
$\bd K$ denote its edge boundary, that is, the set of edges
in $\edges(G)$ having one vertex in $K$ and one outside of
$K$. 
A transitive graph $G$ is {\bf amenable} if 
$\inf |\bd K|/|K| =0$, where $K$ runs over all finite nonempty
vertex sets $K\subset\verts(G)$.

An infinite set of vertices $V_0\subset\verts(G)$
is {\bf end convergent} if for every finite $K\subset\verts(G)$,
there is a component of $G\setminu K$ that contains all
but finitely many vertices of $V_0$.
Two end-convergent sets $V_0,V_1$ are {\bf equivalent}
if $V_0\cup V_1$ is end convergent.
An {\bf end} of $G$ is an equivalence class of end-convergent sets.
Let $\xi$ be an end of $G$.
A {\bf neighborhood} of $\xi$ is a set of vertices in
$G$ that intersects every end-convergent set in $\xi$.
In particular, when $K\subset\verts(G)$ is finite, there is a component
of $G\setminu K$ that is a neighborhood of $\xi$.

Given a set $A$, let $2^A$ be
the collection of all subsets $\eta\subset A$, equipped with the
$\sigma$-field generated by the 
events $\{a\in\eta\}$, where $a\in A$.
A {\bf bond percolation process} is a pair $(\P,\omega)$,
where $\omega$ is a random element in $2^\edges$ and $\P$ denotes
the distribution (law) of $\omega$.  
Sometimes, for brevity, we will just say that $\omega$ is a (bond) percolation.
A {\bf site percolation process} $(\P,\omega)$ is given by a probability
measure $\P$ on $2^{\verts(G)}$, while a (mixed) {\bf percolation}
is given by a probability measure on $2^{\verts(G)\cup\edges(G)}$
that is supported on subgraphs of $G$.
If $\omega$ is a bond percolation process,
then $\hat\omega :=\verts(G)\cup\omega$
is the associated mixed percolation.
In this case, we shall often not distinguish between $\omega$ and
$\hat\omega$, and think of $\omega$ as a subgraph of $G$.
Similarly, if $\omega$ is a site percolation, there is an
associated mixed percolation
$\hat\omega:=\omega\cup\bigl(\edges(G)\cap(\omega\times\omega)\bigr)$,
and we shall often not bother to distinguish between $\omega$ and
$\hat\omega$.

If $v\in\verts(G)$ and $\omega$ is a percolation on $G$,
the {\bf cluster} (or {\bf component}) $C(v)$ of $v$ in $\omega$
is the set of vertices in $\verts(G)$ that can be connected
to $v$ by paths contained in $\omega$.
We will often not distinguish between the cluster $C(v)$ and
the graph
$\Bigl(C(v),\bigl(C(v)\times C(v)\bigr)\cap\omega\Bigl)$ whose vertices
are $C(v)$ and whose edges are the edges in $\omega$ with endpoints
in $C(v)$.

Let $p\in[0,1]$.  Then the distribution of {\bf Bernoulli($p$) bond percolation}
$\omega_p$ on $G$ is the product measure on $2^{\edges(G)}$
that satisfies $\P[e\in\omega]=p$ for all $e\in\edges(G)$.
Similarly, one defines {\bf Bernoulli($p$) site percolation} on
$2^{\verts(G)}$.
The {\bf critical probability} $p_c(G)$ is the infimum over all $p\in[0,1]$
such that there is positive probability for the existence of an infinite
connected component in $\omega_p$.

Aizenman, Kesten and Newman (1987) showed that Bernoulli($p$) percolation
in $\Z^d$ has a.s.\ at most one infinite cluster.
Burton and Keane (1989) gave a much simpler argument that generalizes from
$\Z^d$ to any amenable Cayley graph, though this generalization was not
mentioned explicitly until Gandolfi, Keane and Newman (1992).
It follows that $p_u(G) = p_c(G)$ when $G$ is an amenable Cayley graph.
For background on percolation, especially in $\Z^d$, see Grimmett (1989).

Suppose that $\Gamma$ is an automorphism group of a graph $G$.
A percolation process $(\P,\omega)$ on $G$ is {\bf $\Gamma$-invariant}
if $\P$ is invariant under each $\gamma\in\Gamma$.
This is, of course, the case for Bernoulli percolation.

\bsection {Cluster Indistinguishability}{s.proof} %sz

\procl d.inds %\procname{Indistinguishable clusters}
Let $G$ be graph and $\Gamma$ a closed \hbox{(vertex-)} transitive subgroup
of $\Aut(G)$.  Let $(\P,\omega)$ be a $\Gamma$-invariant bond percolation 
process on $G$.
We say that $\P$ has {\bf indistinguishable infinite clusters}
if for every measurable
$\ev A\subset 2^{\verts(G)}\times 2^{\edges(G)}$
that is invariant under the diagonal action of $\Gamma$,
almost surely, for all infinite clusters $C$ of $\omega$,
we have $(C,\omega)\in\ev A$, or for all
infinite clusters $C$, we have $(C,\omega)\notin\ev A$.
\endprocl

\procl d.ins %\procname{Insertion tolerance}
Let $G=\big(\verts(G),\edges(G)\big)$ be a graph.
Given a set $A\in 2^{\edges(G)}$ and an edge
$e\in\edges(G)$, denote $\ins_e A:=A\cup\{e\}$.
For $\A\subset 2^{\edges(G)}$, we write $\ins_e\ev A := \{\ins_e A \st A \in
\A\}$.
A bond percolation process $(\P,\omega)$ on $G$ is {\bf insertion tolerant} if  
$\P[\ins_e\ev A]>0$
for every $e\in\edges(G)$ and every measurable $\A\subset 2^{\edges(G)}$
satisfying $\P[\A]>0$.
\endprocl

For example, Bernoulli($p$) bond percolation
is insertion tolerant when $p\in \OC{0, 1}$.

A percolation $\omega$  is {\bf deletion tolerant}
if 
$\P[\ins_{\neg e}\ev A]>0$ whenever $e\in\edges(G)$ and $\P[\ev A]>0$,
where $\ins_{\neg e}\omega:=\omega\setminu\{e\}$.
It turns out that deletion tolerance and insertion tolerance
have very different implications.  
For indistinguishability of infinite clusters, we will need insertion
tolerance; 
deletion tolerance does not imply indistinguishability of infinite clusters
(see \ref x.delet/ below).
A percolation that is both insertion and deletion tolerant is usually said
to have ``finite energy''.
Gandolfi, Keane and Newman (1992) use the words ``positive finite energy''
in place of ``insertion tolerance".

\procl t.cerg \procname{Cluster Indistinguishability} 
\assumptions.
Every $\Gamma$-invariant, insertion-tolerant, bond percolation process on $G$
has indistinguishable infinite clusters.
\endprocl

Similar statements hold for site and mixed percolations
and the proofs go along the same lines.  Likewise, the proof extends to
quasi-transitive unimodular automorphism groups.

In a previous draft of our results,
we could only establish the theorem under the assumption
of {\bf strong insertion-tolerance}; that is,
$\P[\ins_e\ev A]\geq\delta\P[\ev A]$ for some constant $\delta>0$.
We are grateful to Olle H\"aggstr\"om for pointing out how to deal with
the general case.

\procl r.cergscen \procname{Scenery}
For some purposes, the following more general form of this theorem
is useful.  Let $G$ be a graph and $\Gamma$
a transitive group acting on $G$. Suppose that $X$ is either $\verts$,
$\edges$, or $\verts\cup\edges$. Let $Q$ be a measurable space and 
$\Omega:=2^\edges\times Q^X$.
A probability measure $\P$ on $\Omega$ will be called a {\bf bond percolation
with scenery} on $G$.  The projection onto $2^\edges$ is the
underlying percolation and the projection onto $Q^{X}$
is the scenery.  If $(\omega,q)\in\Omega$, we set
$\ins_e(\omega,q):= (\ins_e\omega,q)$.  We say that the percolation
with scenery $\P$ is {\bf insertion-tolerant} if, as before,
$\P[\ins_e\ev B]>0$
for every measurable $\ev B\subset\Omega$ of positive measure.
We say that $\P$ has {\bf indistinguishable infinite clusters}
if for every $\ev A \subset 2^\verts\times 2^\edges\times Q^X$
that is invariant under the diagonal action of $\Gamma$, for $\P$-a.e.\
$(\omega, q)$, either
all infinite clusters $C$ of $\omega$ satisfy
$(C,\omega, q)\in\ev A$ or they all satisfy $(C,\omega, q)\notin\ev A$.
\ref t.cerg/ holds also for percolation with scenery
(with the same proof as given below).

A percolation with scenery that comes up naturally is as follows.
Fix $p_1<p_2$ in $(0,1)$.  Let $z_e$ ($e\in\edges$) be i.i.d.\ 
with each $z_e$ distributed according to uniform measure on $[0,1]$.
Let $\omega_j$ be the set of $e\in\edges$ with $z_e<p_j$, $j=1,2$.
Then $(\omega_2,\omega_1)$ is an insertion-tolerant percolation with
scenery, where $\omega_1\in 2^\edges$ is considered the scenery.
See H\"aggstr\"om and Peres (1999),
H\"aggstr\"om, Peres and Schonmann (1998),
Schonmann (1999a), and Alexander (1995) for examples where this process is
studied.
\endprocl

Say that an infinite cluster $C$ is of {\bf type $\ev A$} if
$(C,\omega)\in\ev A$;
otherwise, say that it is of {\bf type $\neg\ev A$}.
Suppose that there is an infinite cluster
$C$ of $\omega$ and an edge $e\in\edges$
with $e\notin \omega$ such that the
connected component $C'$ of $\omega\cup\{e\}$
that contains $C$ has a type different from the type of $C$.
Then $e$ is called {\bf pivotal} for $(C,\omega)$.

We begin with an outline of the proof of \ref t.cerg/.
Assume that the theorem fails.  First, we shall show that with positive
probability, given that a vertex belongs to an infinite cluster,
there are pivotal edges at some distance $r$ of the vertex. 
By \ref p.trans/, a.s.\ when there is more
than one infinite cluster of $\omega$, each such cluster is transient
for the simple random walk, and hence for the so-called delayed simple random
walk (DSRW).  The DSRW on a cluster of $\omega$ is stationary (in the
sense of \ref l.dsrw/).
Fix a base point $o\in\verts$, and let $W$ be the
DSRW on the $\omega$-cluster of $o$ with $W(0)=o$.
Let $n$ be large, and let $e$ be a uniform random edge at distance $r$
{}from $W(n)$.  When DSRW is transient, with probability bounded away from
zero, $e$ is pivotal for $\big(C(o),\omega\big)$
and $W(j)$ is not an endpoint of $e$ for
any time $j<n$.  On this event, set $\omega':=\ins_e\omega$.
Then $\omega'$ is, up to a controllable factor, as likely as $\omega$
by insertion tolerance.  By transience, $e$ is far from $o$ 
with high probability.
Since $e$ is pivotal, the type of $C(o)$ is different in $\omega$ and
in $\omega'$.  Since $\omega$ and $\omega'$ are the same
in a large neighborhood of $o$, this shows that the type
of $C(o)$ cannot be determined with arbitrary accuracy
by looking at $\omega$ in a large neighborhood of $o$.
This contradicts the measurability of $\ev A$ and
establishes the theorem.

One can say that the proof is based on the
contradictory prevalence of pivotal edges.  To put the situation in the correct
perspective, we point out that there are events depending on i.i.d.~zero-one
variables that have infinitely many ``pivotals'' with positive probability.
For example, let $m_1,m_2,\dots$ be a sequence of positive integers such that
$\sum_k 2^{-m_k}<\infty$ but $\sum_k m_k 2^{-m_k}=\infty$, and
let $\Seq{x_{k,j}\st j\in\{1,\dots,m_k\}}$ be i.i.d.\ 
random variables with $\P[x_{k,j}=0]=\P[x_{k,j}=1]=1/2$.  Let
$\ev X$ be the event that there is a $k$ such that
$x_{k,j}=1$ for $j=1,2,\dots,m_k$.  Then with positive probability
$\ev X$ has infinitely many pivotals, i.e., there are infinitely many $(k,
j)$ such that changing $x_{k, j}$ from 0 to 1 will change from $\neg
\ev X$ to $\ev X$.  (This is a minor variation on an example
described by H\"aggstr\"om and Peres (private communication).)

We now prove the lemmas necessary for \ref t.cerg/.

\procl l.piv \procname{Pivotals}
Suppose that $(\P,\omega)$ is an insertion-tolerant percolation
process on a graph $G$, and let $\ev A\subset 2^{\verts(G)}\times 2^{\edges(G)}$
be measurable.  Assume that there is positive probability for
co-existence of infinite clusters in $\ev A$ and $\neg \ev A$.
Then with positive probability, there is an infinite cluster $C$ of the percolation
that has a pivotal edge.
\endprocl

\proof
Let $k$ be the least integer such that
there is positive probability that there
are infinite clusters of different types with distance
between them equal to $k$.
Clearly, $k>0$.  
Suppose that $\gamma$ is a path of length $k$
such that with positive probability, $\gamma$
connects infinite clusters of different types;
let $\ev G$ be the event that
$\gamma$ connects infinite clusters of different types. 
Given $\ev G$, there
are exactly two infinite clusters that intersect $\gamma$, by the
minimality of $k$.
Let $e$ be the first edge in $\gamma$.
When $\omega\in\ev G$, let $C_1$ and $C_2$
be the two infinite clusters that $\gamma$
connects, and let $C_1'$ and $C_2'$ be the
infinite clusters of $\ins_{e}\omega$ that contain
$C_1$ and $C_2$, respectively.
Note that, conditioned on $\ev G$,
the distance between $C_1'$ and $C_2'$ is less than $k$.
(The possibility that $C_1'=C_2'$ is allowed.)
Since $\ins_{e}\ev G$ has positive probability,
the definition of $k$ ensures that the type of $(C_1',\ins_e\omega)$
equals the type of $(C_2',\ins_e\omega)$ a.s.
Hence, $e$ is pivotal for $(C_1,\omega)$ or for
$(C_2,\omega)$ whenever $\omega\in\ev G$.
This proves the lemma.
\Qed

\procl l.ergcomp
\genassumptions.
If $(\P,\omega)$ is a $\Gamma$-invariant insertion-tolerant percolation
process on an infinite graph $G$, then almost every ergodic component of $\omega$
is insertion tolerant.
\endprocl

\procl r.GKN A stronger and more general statement is Lemma 1 of Gandolfi,
Keane and Newman (1992).
\endprocl

\proof 
This is the same as saying that for every $\Gamma$-invariant event $\ev G$
of positive probability, the probability measure $\P[\;\cbuldot \mid \ev G]$
is insertion tolerant.
If $\ev G$ is a tail event, then $\ins_e\ev G=\ev G$, and insertion tolerance
of $\P[\;\cbuldot \mid \ev G]$ follows from the calculation
$$
\P[\ins_e \ev A \mid \ev G]
=
{\P[\ins_e \ev A \cap \ev G] \over \P[\ev G]}
=
{\P[\ins_e \ev A \cap \ins_e \ev G] \over \P[\ev G]}
=
{\P[\ins_e (\ev A \cap \ev G)] \over \P[\ev G]}
\,.  
$$
But every $\Gamma$-invariant event is a tail event (mod 0): Write $\edges$ as an
increasing union of finite subsets $\edges_n$. Let $\ev G_n$ be events that
do not depend on edges in $\edges_n$ and such that $\sum_n
\P[\ev G_n \xor \ev G] < \infty$; such events exist because $\P$ and
$\ev G$ are $\Gamma$-invariant.
(Here $\ev H\xor\ev K:= (\ev H\setminus\ev K)\cup(\ev K\setminus\ev H)$
is the exclusive or.)
Then $\limsup_n \ev G_n$ is a tail event and equals $\ev G$ (mod 0). \Qed

\procl c.ninf % \procname{Number of infinite clusters}
\genassumptions.
Consider a $\Gamma$-invariant insertion-tolerant percolation
process $(\P,\omega)$ on $G$. Then almost surely, the number of infinite
clusters of $\omega$ is $0,1$ or $\infty$.
\endprocl

The proof is standard for the ergodic components; cf.\ Newman and
Schulman (1981). 

Benjamini and Schramm (1996) conjectured that for Bernoulli percolation on
any quasi-transitive graph, if there are infinitely many infinite clusters,
then a.s.\ every infinite cluster has continuum many ends. This was proved
by H\"aggstr\"om and Peres (1999) in the unimodular case and then by
H\"aggstr\"om, Peres, and Schonmann (1998) in general.
For the unimodular case, we give a simpler proof that extends to
insertion-tolerant percolation processes.
We begin with the following proposition.

\procl p.niso
\assumptions.
Consider some $\Gamma$-invariant percolation process $(\P,\omega)$ on $G$.
Then a.s.\ each infinite cluster that has at least $3$ ends has
no isolated ends.
\endprocl

\proof
For each $n=1,2,\dots$, let $\hat A_n$ be the union of
all vertex sets $A$ that are contained in some percolation cluster $K(A)$,
have diameter at most $n$ in the metric of the percolation cluster,
and such that $K(A)\setminu A$ has at least $3$ infinite components.
Note that if $\xi$ is an isolated end of a percolation cluster $K$,
then for each finite $n$, some neighborhood of $\xi$ in $K$ is
disjoint from $\hat A_n$.
Also observe that
if $K$ is a cluster with at least 3 ends, then
$K$ intersects $\hat A_n$ for some $n$.

Fix some $n \ge 1$.
Consider the mass transport that sends one unit of mass from each vertex $v$
in a percolation cluster that intersects $\hat A_n$ and distributes 
it equally among the
vertices in $\hat A_n$ that are closest to $v$ in the metric of the percolation
cluster of $v$. 
In other words, let $K(v)$ be the set of vertices
in $C(v)\cap\hat A_n$ that are closest to $v$ in the metric of $\omega$,
and set $F(v, w; \omega):= \left| K(v)\right|^{-1}$ if $w\in K(v)$ and
otherwise $F(v,w;\omega):= 0$.
Then $F(v, w; \omega)$, and hence the expected mass $f(v, w) := \E F(v, w;
\omega)$ transported from $v$ to $w$, is invariant under the diagonal
$\Gamma$ action.
If $\xi$ is an isolated end of an infinite cluster $K$ that intersects
$\hat A_n$, then there is a finite set of vertices $B$ that
gets the mass from all the vertices in a neighborhood of $\xi$.
But the Mass-Transport Principle tells us that the expected
mass transported to a vertex is finite.  Hence, a.s.\ clusters that
intersect $\hat A_n$ do not have isolated ends.  Since this holds
for all $n$, we gather that a.s.\ infinite clusters with isolated ends
do not intersect $\bigcup_n\hat A_n$, whence they have at most two ends.
 \Qed

\procl p.isolated
\assumptions.
If $(\P,\omega)$ is a $\Gamma$-invariant insertion-tolerant percolation process on
$G$ with infinitely many infinite clusters a.s., then a.s.\ every infinite
cluster has continuum many ends and no isolated end.
\endprocl

\proof
As Benjamini and Schramm (1996) noted, it suffices to prove that there are
no isolated ends of clusters.
To prove this in turn, observe that if some cluster has an isolated end, then
because of insertion tolerance, with positive probability, some cluster
will have at least 3 ends with one of them being isolated.  Hence \ref
p.isolated/ follows from \ref p.niso/.
\Qed

\procl p.trans \procname{Transience}
\assumptions.
Suppose that $(\P,\omega)$ is a $\Gamma$-invariant insertion-tolerant
percolation process on $G$ that has almost surely infinitely many
infinite clusters.  Then a.s.\ each infinite cluster is transient.
\endprocl

\proof
By \ref p.isolated/,
every infinite cluster of $\omega$ has infinitely many ends. 
Consequently, there is a random forest $\fo\subset\omega$ whose distribution
is $\Gamma$-invariant such that a.s.\ each infinite cluster $C$ of
$\omega$ contains a tree of $\fo$ with more than $2$ ends (Lemma 7.4 of
\BLPSgip).
{}From Remark 7.3 of \BLPSgip, we know that any such tree has $p_c<1$.
By Lyons (1990), it follows that such a tree is transient. 
The Rayleigh monotonicity principle
(e.g., Lyons and Peres (1998)) then implies that $C$ is transient.
\Qed

\procl r.trns-nonuni
Examples show that \ref p.trans/ does not hold when $\Gamma$ is not
unimodular. However, we believe that when $\Gamma$ is not unimodular and
$\P$ is Bernoulli percolation, it does still hold.
\endprocl

Let $(\P,\omega)$ be a bond percolation process on $G$, and let $\omega\in2^\edges$.
Let $x\in\verts(G)$ be some base point.
It will be useful to consider {\bf delayed simple random walk} on $\omega$
starting at $x$, $W = W_x^\omega$, defined as follows.
Set $W(0):=x$.
If $n\geq 0$, conditioned on $\Seq{W(0),\dots,W(n)}$ and $\omega$,
let $W'(n+1)$ be chosen from the neighbors of $W(n)$ 
in $G$ with equal probability.
Set $W(n+1) := W'(n+1)$ if the edge $[W(n), W'(n+1)]$ belongs to $\omega$;
otherwise, let $W(n+1) := W(n)$.

Given $\omega$, let $W$ and $W^*$ be two independent delayed simple random
walks starting at $x$.  Set $w(n):= W(n)$ for $n\geq 0$ and $w(n) := W^*(-n)$
for $n<0$. Then $w$ is called {\bf two-sided} delayed simple random walk.
Let $\Ph_x$ denote the law of the pair $(w,\omega)$; it
is a probability measure on $\verts^\Z\times2^\edges$.
Define $\SS:\verts^\Z\to\verts^\Z$ by
$$
(\SS w)(n):= w(n+1)
\,,
$$
and let 
$$
\SS(w,\omega):= (\SS w,\omega) \qquad\forall (w,\omega)\in\verts^\Z\times2^\edges
\,.
$$
For $\gamma\in\Gamma$, we set
$$
\gamma(w,\omega):= (\gamma w,\gamma\omega)
\,,
$$ 
where $(\gamma w)(n):=\gamma\bigl(w(n)\bigr)$.

The following lemma generalizes similar lemmas in H\"aggstr\"om (1997),
in H\"aggstr\"om and Peres (1999), and in Lyons and Peres (1998).  

\procl l.dsrw
\procname{Stationarity of Delayed Random Walk}
\assumptions.
Let $o\in\verts(G)$ be some base point.
Let $(\P,\omega)$ be a $\Gamma$-invariant bond percolation process on $G$.
Let $\Ph_o$ be the joint law of $\omega$ and two-sided delayed simple random walk
on $\omega$, as defined above.  
Then $\Ph_o[\ev A] = \Ph_o[\SS\ev A]$ for every $\Gamma$-invariant
$\ev A\subset \verts^\Z\times 2^\edges$.
In other words, the restriction of $\Ph_o$ to the $\Gamma$-invariant $\sigma$-field
is $\SS$-stationary.
\endprocl

The lemma will follow from two identities.  The first is based on the
fact that the transition operator for delayed simple random walk on $\omega$
is symmetric,
and the second is based on the Mass-Transport Principle, i.e., on
unimodularity.

\proof
For $j\in\Z$ and $x\in\verts$, set
$$
\ev W^j_x:= \left\{(w,\omega)\in\verts^\Z\times 2^\edges\st w(j)=x\right\}
$$
and
$$
\mu:= \sum_{x\in\verts}\Ph_x
\,.
$$
Let $\ev B\subset\verts^\Z\times2^\edges$ be measurable.
Observe that for all $\omega\in2^\edges$,
$$
\mu[\ev B \mid \omega] = \mu[\SS\ev B\mid\omega]
\,,
\label e.ssinv
$$
where $\mu[\ev B \mid \omega]$ means $\sum_{x\in\verts}\Ph_x[\ev B\mid\omega]$.
This follows from the fact that for any $j,k\in\Z$ with $j<0<k$
and any $v_{j},v_{j+1},\dots,v_k\in\verts$,
we have 
$$
\mu\left[
\ev W^{j}_{v_j}\cap \ev W^{j+1}_{v_{j+1}}\cap \cdots \cap \ev W^k_{v_k}\mid \omega\right]
= \prod_{i=j}^{k-1} a_i
\,,
$$
where $a_i:=\deg_G(v_i)^{-1}=\deg_G(o)^{-1}$ if $[v_i,v_{i+1}]\in\omega$, 
$a_i := \bigl(\deg_G(o)-\deg_\omega(v_i)\bigr)/\deg_G(o)$ if
$v_i=v_{i+1}$, and $a_i:=0$ otherwise.
By integrating \ref e.ssinv/ over $\omega$, we obtain
$$
\mu [\ev B] = \mu[\SS\ev B]
\,,
\label e.fiden
$$
which is our first identity.  Observe also that 
$\mu$ is $\Gamma$-invariant.

Now let $\ev A\subset\verts^\Z\times2^\edges$ be $\Gamma$-invariant and
measurable.
For $x,y\in\verts$, define 
$$
\phi(x,y) := \mu[\ev A\cap \ev W^{-1}_x\cap \ev W^0_y]
\,.
$$
Then $\phi$ is invariant under the diagonal
action of $\Gamma$ on $\verts\times\verts$,
because $\mu$ and $\ev A$ are $\Gamma$-invariant.
Consequently, the Mass-Transport Principle gives
$\sum_{x\in\verts}\phi(x,o)=\sum_{x\in\verts}\phi(o,x)$,
which translates to our second identity:
$$
\mu[\ev A \cap \ev W^0_o] = \mu[\ev A\cap\ev W^{-1}_o]
\,.
\label e.siden
$$

Observe that $\mu[\ev C\cap\ev W^0_o] = \Ph_o[\ev C]$ for all
measurable $\ev C\subset\verts^\Z\times2^\edges$.
By using \ref e.fiden/ with $\ev B:= \ev A\cap\ev W^0_o$ and
then using \ref e.siden/ with $\SS\ev A$ in place of $\ev A$,
we obtain finally
$$
\Ph_o[\ev A] =
\mu[\ev A \cap \ev W^0_o] = \mu\bigl[ \SS (\ev A \cap\ev W^0_o)\bigr]
=
\mu\bigl[\SS \ev A \cap \ev W^{-1}_o\bigr]
=
\mu\bigl[\SS \ev A \cap \ev W^{0}_o\bigr]
= 
\Ph_o[\SS\ev A]
\,.
\Qed
$$ %sz

In \ref t.ergdel/ below, we show that in an appropriate sense,
if $\omega$ is ergodic and has indistinguishable components, then the
delayed simple random walk on the infinite components of $\omega$ is
ergodic.

\procl r.dsrwgen
\procname{A generalization of \ref l.dsrw/}
Let $\verts$ be a countable set acted on by a transitive unimodular group
$\Gamma$.
Let $Q$ be a measurable space, 
let $\Xi := Q^{\verts\times\verts}$, and
let $\P$ be a $\Gamma$-invariant probability measure on $\Xi$.
Suppose that $z:Q\to[0,1]$ is measurable, and that
$z\bigl(\xi(x,y)\bigr)=z\bigl(\xi(y,x)\bigr)$ and
$\sum_v z\bigl(\xi(x,v)\bigr)=1$ for all $x,y\in\verts$ and for
$\P$-a.e.~$\xi$.
Given $o\in\verts$ and a.e.~$\xi\in\Xi$,
there is an associated random walk starting at $o$ with transition
probabilities $p_\xi(x,y)=z\bigl(\xi(x,y)\bigr)$.  Let $\Ph_o$ denote the joint
distribution of $\xi$ and this random walk.  Then the above
proof shows that the restriction
of $\Ph_o$ to the $\Gamma$-invariant $\sigma$-field is $\SS$-invariant.
\par
See Lyons and Schramm (1998) for a still greater generalization.
\endprocl

Let $(\P,\omega)$ be a bond percolation process on an infinite graph $G$.
For every $e\in\edges(G)$, let $\ev F_{\neg e}$ be the $\sigma$-field
generated by the events $\{e'\in\omega\}$ with $e'\neq e$.
Set $Z(e):=\P[e\in\omega\mid \ev F_{\neg e}]$,
and call this the {\bf conditional marginal} of $e$.
Note that insertion tolerance is equivalent to $Z(e)>0$ a.s.\ for
every $e\in\edges(G)$.  

\proofof t.cerg 
Let $(\P,\omega)$ be insertion tolerant and $\Gamma$-invariant.
Let $o\in\verts$ be some fixed base point of $G$.
Assume that the theorem is false.  Then by \ref c.ninf/,
there is positive probability that $\omega$ has
infinitely many infinite clusters.  If we
condition on this event, then $\omega$ is still insertion tolerant, as
shown in \ref l.ergcomp/. 
Consequently, we henceforth assume, with no loss of generality,
that a.s.\ $\omega$ has infinitely many infinite clusters.

Fix $\epsilon>0$. 
{}From \ref l.piv/, we know that there is a positive probability
for pivotal edges of clusters of type $\ev A$, or
there is a positive probability for pivotal edges of clusters 
of type $\neg\ev A$.  Since we may replace $\ev A$ by its complement,
assume, with no loss of generality, that there is
a positive probability for pivotal edges of clusters of type
$\ev A$.
Fix some $r>0$ and $\delta>0$ such that with positive probability,
$C(o)$ is infinite of type $\ev A$ and
there is an edge $e$ at distance $r$ from $o$
that is pivotal for $\big(C(o),\omega\big)$ and
satisfies $Z(e)>\delta$, where $Z(e)$ is the conditional marginal of $e$,
as described above.

Let $\ev A_o$ be the event that $C(o)$ is infinite and
$\big(C(o),\omega\big)\in\ev A$,
and let $\wtA $ be an event that depends on only finitely
many edges such that $\P[\ev A_o\xor \wtA ]<\epsilon$.
Let $R$ be large enough that $\wtA $ depends only
on edges in the ball $B(o,R)$.

Let $W:\Z \to\verts$
be two-sided delayed simple random walk on $\omega$, with $W(0)=o$.
For $n\in\Z$, let $e_n\in\edges$ be an edge chosen uniformly among the
edges at distance $r$ from $W(n)$.
Write $\Ph$ for the probability measure where we choose $\omega$ according
to $\P$, choose $W$, and choose $\Seq{e_n\st n\in\Z}$ as indicated.

Given any $e\in\edges$, 
let $\ev P_{e}$ be the event that $\omega\in\ev A_o$, that
$e$ is pivotal for $C(o)$, and that $Z(e)>\delta$.
Let $\ev E^n_{e}$ be the event that $e_n=e$ and $W(j)$ is not an endpoint of
$e$ whenever $-\infty<j<n$.
Note that for all $\zeta\in 2^\edges$, $n \in \Z$, and $e\in\edges$,
$$
\Ph[\ev E^n_{e} \mid \omega=\zeta] =
\Ph[\ev E^n_{e} \mid \omega=\ins_{e}\zeta]
\,.
$$
Thus, for all measurable $\ev B\subset2^\edges$ with $\P[\ev B]>0$
and $\P[Z(e)\geq\delta\mid\ev B]=1$, we have
$$
\Ph[\ev E^n_{e} \cap \ins_e \ev B] =
\Ph[\ev E^n_{e}  \mid \ins_e \ev B] \P[\ins_e \ev B] =
\Ph[\ev E^n_{e}  \mid  \ev B] \P[\ins_e \ev B]
% = {\Ph[\ev E^n_{e}  \cap  \ev B] \P[\ins_e \ev B] \over \P[\ev B]}
 = {\P[\ins_e \ev B] \over \P[\ev B]} \Ph[\ev E^n_{e} \cap \ev B]
\ge
\delta \Ph[\ev E^n_{e} \cap \ev B]
\,.
$$
In particular,
$$
\Ph[\ev E^n_{e}\cap\ins_e \wtA  \cap \ins_{e}\ev P_{e}]
\geq \delta \Ph[\ev E^n_{e}\cap\wtA \cap \ev P_{e}]
= \delta \Ph[\ev E^n_{e}\cap\wtA \cap \ev P_{e_n}]
\,.
\label e.aa
$$
Let 
$$
\displaylines{
\ev E^n := \bigcup_{e\in\edges} \ev E^n_{e}
\,,\cr
\ev E_R^n := \bigcup_{e\in\edges\setminu B(o,R)} \ev E^n_{e}
\,, \cr
}$$
and note that these are disjoint unions.
Since $\ins_{e} \ev P_{e}\subset \neg \ev A_o$
and since $\ins_e \wtA \subset \wtA $ when $e\notin B(o,R)$, 
we may sum \ref e.aa/ over all $e\notin B(o,R)$ to obtain that
$$
\P[\wtA  \cap \neg \ev A_o]
\geq \Ph[\ev E_R^n\cap \wtA  \cap \neg \ev A_o]
\geq \delta \Ph[\ev E_R^n \cap \wtA  \cap\ev P_{e_n}]
\geq \delta \Ph[\ev E_R^n \cap \ev A_o \cap\ev P_{e_n}]
-\delta\epsilon
\,.
\label e.bb
$$
Fix $n$ to be sufficiently large that the probability
that $C(o)$ is infinite and $e_n\in B(o,R)$ is smaller than
$\epsilon$;  this can be done by \ref p.trans/.
Then $\Ph[\ev A_o\cap\ev E^n\setminu\ev E_R^n]\leq \epsilon$.
Hence, we have from \ref e.bb/ that
$$
\epsilon\geq \P[\wtA\xor\ev A_o]
\geq \P[\wtA  \cap \neg \ev A_o]
\geq \delta \Ph[\ev E^n\cap\ev A_o\cap\ev P_{e_n}]-
2\epsilon\delta
\,.
\label e.cc
$$ %sz
Recall that $\Ph[\ev A_o\cap\ev P_{e_0}]>0$.
Moreover, conditioned on $\ev A_o\cap\ev P_{e_0}$, transience guarantees
a.s.\ a least $m\in\Z$ such that
$W(m)$ is at distance $r$ to $e_0$.
Consequently, for some $m\leq 0$,
$$
\Ph[\ev E^m\cap\ev A_o\cap\ev P_{e_m}]>0
\,.
\label e.flip
$$
Let $\ev B_m$ be the event that $Z(e_m)\geq\delta$, that
$e_m$ is not an endpoint
of $W(j)$ for $j<m$, that $C\big(W(m)\big)$ is infinite and of type
$\ev A$, and that $e_m$ is pivotal for $C\big(W(m)\big)$.
Then $\ev B_m= \ev E^m\cap\ev A_o\cap\ev P_{e_m}$ up to zero
$\Ph$-measure.  But $\ev B_m$ is $\Gamma$-invariant (up to zero $\Ph$-measure).
Therefore, \ref l.dsrw/ shows that the left-hand side of \ref e.flip/
does not depend on $m$, and it certainly does not depend on $\epsilon$.
Hence, when we take $\epsilon$ to be a sufficiently
small positive number, \ref e.cc/ gives a contradiction.
This completes the proof of the theorem.
\Qed

\procl x.delet
A deletion-tolerant process that does not have
indistinguishable clusters is obtained as follows.
Let $X$ be a 3-regular tree and $p \in (1/2, 1)$.
Let $p' \in (1/(2p), 1)$.
Begin with Bernoulli($p$) percolation on $X$.
Independently for each cluster $C$, with probability 1/2,
intersect it with an independent Bernoulli($p'$) percolation.
The resulting percolation process is clearly deletion tolerant, yet some
infinite clusters $\omega$ have $p_c(\omega) = 1/(2p)$, while others have
$p_c(\omega) = 1/(2pp')$.
\endprocl

\procl r.robust
Let $T$ be the $3$-regular tree, let $\Gamma$ be the subgroup of
automorphisms of $T$ that fixes an end $\xi$ of $T$, and let
$\P$ be Bernoulli($p$) bond percolation on $T$, where $p\in (1/2,1)$.
Then a.s.\ each infinite cluster $C$ has a unique vertex $v_C$ ``closest''
to $\xi$.
The degree of $v_C$ in the percolation configuration distinguishes
among the infinite clusters.  Hence,
\ref t.cerg/ does not hold without the assumption that $\Gamma$ is unimodular.
By a simple modification, a similar example can be constructed
where $\Gamma$ is the full automorphism group of a graph on
which the percolation is performed (compare \BLPSgip).

However, H\"aggstr\"om, Peres and Schonmann (1998) have recently shown
that even without the unimodularity assumption, when $p>p_c$
so-called ``robust'' properties
do not distinguish between the infinite clusters of Bernoulli($p$) percolation.
\endprocl

\procl q.heavy
In the case that $\Gamma$ is nonunimodular, write $\mu_x$ for the Haar
measure of the stabilizer of $x \in \verts$. The infinite clusters $C$ divide
into two types: the {\bf heavy} clusters for which $\sum_{x \in C} \mu_x =
\infty$ and the others, the {\bf light} clusters. It can be that the light
clusters are distinguishable: e.g., consider a 3-regular tree $T$ with a
fixed end $\xi$. Let $\Gamma$ be the group of automorphisms of $T$ that fix
$\xi$. Then Bernoulli(2/3) percolation has infinitely many light clusters,
which can be distinguished by the degree of the vertex they contain that is
closest to $\xi$. But is it the case that heavy clusters are
indistinguishable for every insertion-tolerant percolation process that is
invariant with respect to a transitive automorphism group?
\endprocl

\bsection {Uniqueness and Connectivity}{s.decay}

Our goal is to prove

\procl t.tau
\procname{Uniqueness and Connectivity}
\assumptionsandinf.
Let $\P$ be a $\Gamma$-invariant and ergodic insertion-tolerant
percolation process on $G$. If $\P$ has more than one infinite
component a.s., then connectivity decays:
$$
\inf\bigl\{ \tau(x, y)\st{x, y \in \verts} \bigr\} = 0
\,.
$$
\endprocl

The intuitive idea behind our proof is that if $\inf \tau(x, y) > 0$, then
each infinite cluster has a positive ``density''.
Since the densities are the same by cluster indistinguishability, there are
only finitely many infinite clusters.
By \ref c.ninf/, there is only one.
To make the idea of ``density'' precise, we use simple random walk $X$ on the
whole of $G$ with the percolation subgraph as the scenery, counting how
many times we visit each cluster.

For a set $C\subset\verts$, write 
$$
\freq(C) := \lim_{n \to\infty} {1 \over n} \sum_{k = 1}^n \II{X(k) \in C}
\,,
$$ %sz
when the limit exists,
for the frequency of visits to $C$ by the simple random walk on $G$.

\procl l.freq \procname{Cluster Frequencies}
\assumptions.
There is a $\Gamma$-invariant measurable function $\freqfn:2^\verts\to[0,1]$
with the following property.
Suppose that $(\P,\omega)$ is a $\Gamma$-invariant bond percolation process on $G$,
and let $\Ph:=\P\times\P_o$, where $\P_o$ is the law of simple random
walk on $G$ starting at the base point $o$.
Then $\Ph$-a.s., $\freq(C)=\freqfn(C)$ for every cluster $C$.
\endprocl

\proof
Given a set $C\subset\verts$ and $m,n\in\Z$, $m< n$, let
$$
\freq^n_m(C) := {1 \over n-m} \sum_{k=m}^{n-1} \II{X(k) \in C}
\,.
$$
For every $\freq\in[0,1]$, let 
$$
\ev Z_\freq:= \left\{C\subset\verts\st
\lim_{n\to\infty} \freq^n_0(C)=\alpha \ \P_o\hbox{-a.s.{}}\right\}
\,,
$$ 
and set $\ev Z := \cup_{\freq\in[0,1]} \ev Z_\freq$.
Define $\freqfn(C):=\freq$ when $C\in\ev Z_\freq$.
If $C\notin\ev Z$, put $\freqfn(C):=0$, say.
It is easy to verify that $\freqfn$ is measurable.
Let $\gamma\in\Gamma$.  To prove that a.s.\ $\freqfn(C)=\freqfn(\gamma C)$,
note that there is an $m\in\N$ such that with positive
probability $X(m)=\gamma o$.  Hence for every
measurable $A\subset [0,1]$ such that $\freq(C)\in A$ with
positive probability, we have $\freq (\gamma C)\in A$ with positive
probability.  This implies that $\freqfn$ is $\Gamma$-invariant.

It remains to prove that $\P$-a.s., every component of $\omega$ is
in $\ev Z$.
First observe that the restriction of $\Ph$ to the 
$\Gamma$-invariant $\sigma$-field is $\SS$-invariant.
This can be verified directly, but is also special case of \ref r.dsrwgen/:
take $Q=\{-1,0,1\}$, take the value in $Q$ associated to a
pair $(x,y)\in\verts\times\verts$ to be $1$ if $[x,y]\in\omega$, $0$ if
$[x,y]\in\edges\setminu\omega$, and $-1$ if $[x,y]\notin\edges$, and take
$z(0)=z(1)=1/\deg_G(o)$ and $z(-1)=0$.

The following argument is modeled on the proof of Thm.~1 of
Burton and Keane (1989).
Let $F^n(j)$ be
the number of times that the $j$ most frequently visited clusters
in $[0,n-1]$ are visited in $[0,n-1]$.
That is, 
$$
F^n(j) := n\, \max\bigl\{
\freq^n_0(C_1)+\cdots+\freq^n_0(C_j)
\st C_1, \ldots, C_j \hbox{ are distinct clusters} \bigr\}\,.
$$  %sz
For each fixed $j$, it is easy to see that $F^n(j)$ is a subadditive
sequence, i.e., $F^{n+k}(j) \le F^n(j) + \SS^n F^k(j)$.
Note that the random variables $F^n(j)$ are
invariant with respect to the diagonal action of $\Gamma$
on $2^\edges \times \verts^\Z$. 
By the subadditive ergodic theorem, $\lim_n F^n(j)/n$
exists a.s.
Set $\freq(j):=\lim_n F^n(j)/n-\lim_n F^n(j-1)/n$ for $j \ge 1$.
We claim that a.s.
$$
\lim_{n\to\infty} \max \bigl\{ |\freq^m_0(C)-\freq^k_0(C)|
\st k,m\in\Z,\,k,m\geq n,\, C\hbox{ is a cluster}\bigr\} = 0
\,,
\label e.ucon
$$ %sz
the Cauchy property uniform in $C$.
Indeed, let $\eps>0$.  Observe that $\sum_j\alpha(j)\leq 1$
and that $\alpha(1)\geq\alpha(2)\geq\cdots$.
Let $j_1\geq 1$ be large enough that $\alpha(j)<\eps/9$ for all $j\geq j_1$.
Let $m_1$ be sufficiently large that
$\bigl|F^m(j)/m-F^m(j-1)/m-\alpha(j)\bigr|<\eps/(9j_1+9)$
for all $j\leq j_1$ and all $m\geq m_1$.
Set 
$$
U:=\bigl\{x\in[0,1]\st \texists{j\leq j_1}\left|x-\alpha(j)\right|
<\eps/(9 j_1+9)\bigr\} \cup [0,\eps/3]
\,,
$$
and let $U(\delta)$ denote the set of points $x\in\R$ within
distance $\delta$ of $U$.
For all clusters $C$ that are visited by the random walk and all $m=1,2,\dots$,
there is some $j$ such that $\alpha_0^m(C)=F^m(j)/m-F^m(j-1)/m$.
If $j\leq j_1$ and $m\geq m_1$, it follows that $\alpha_0^m(C)\in U$.
The same also holds when $j\geq j_1$ and $m\geq m_1$, because the $j_1$-th most
frequently visited cluster in $[0,m-1]$, say $C'$,
satisfies $\alpha_0^m(C)\leq\alpha_0^m(C')\leq \alpha(j_1)+\eps/9\leq \eps/3$.
Hence $\alpha_0^m(C)\in U$ for all clusters $C$ and all $m\geq m_1$.
Because $-1/m\leq\alpha_0^m(C)-\alpha_0^{m+1}(C)\leq1/m$,
it follows that for all $C$ and all $n\geq m_1$, the set
$\{\alpha_0^m(C)\st m\geq n\}$ is contained in
some connected component of  $U({1/n})$.
But when $n>\max\{m_1,9(j_1+1)/\eps\}$, the total length of
$U(1/n)$ is less than $\eps$,
which implies that the diameter of each connected component
is less than $\eps$.
This verifies \ref e.ucon/.

Since 
$$
2\max \bigl\{|\freq^{2n}_0(C)-\freq^n_0(C)|\st C\hbox{ is a cluster}\bigr\}
=
\max \bigl\{|\freq^{2n}_n(C)-\freq^n_0(C)|\st C\hbox{ is a cluster}\bigr\}
$$
has the same law as 
$\max\bigl\{|\freq^n_{0}(C) - \freq^0_{-n}(C)|\st C\hbox{ is a cluster}\bigr\}$
(by the $\SS$-invariance noted above),
it follows from  \ref e.ucon/ that a.s.\
$\lim_{n\to\infty} \freq^n_0(C)=\lim_{n\to \infty}\freq^0_{-n}(C)$
for every cluster $C$.
When the cluster $C$ is fixed, $\freq^n_0(C)$ and
$\freq^0_{-n}(C)$ are independent, but both tend to 
$\freq(C)$.  Hence $\freq(C)$ is an a.s.\ constant, 
which means that $\P$-a.s.\  we have $C\in\ev Z$
for every cluster $C$. This completes the proof.
\Qed

\proofof t.tau
Let $\freqfn$ be as in \ref l.freq/.
Since $\P$ is ergodic, \ref t.cerg/ implies that there
is a constant $c\in[0,1]$ such that
a.s.\ $\freqfn(C)=c$ for every infinite cluster $C$.  Suppose that
there is more than one infinite cluster with positive probability. 
Then there are infinitely many infinite clusters a.s.
Since clearly $\sum_C \freq(C)=\sum_C \freqfn(C)\leq 1$,
where the sum is over all clusters, it follows that $c=0$.
Since $G$ is infinite, it is also immediate that
$\freqfn(\{v\})=0$ for every $v\in\verts$.
Therefore, $\freqfn(C)=0$ a.s.\ for all clusters, finite or infinite.
In particular, $\freqfn\big(C(o)\big)=0$ a.s.,
where $C(o)$ is the cluster of $o$.

Let $\tau_0 := \inf\bigl\{ \tau(x, y)\st x,y\in\verts\bigr\}$. 
We have that
$$
\E[\freq_0^n\big(C(o)\big)] \ge \tau_0\,,
$$
whence 
$0=\E[\freqfn\big(C(o)\big)]= \E[\freq\big(C(o)\big)] \ge \tau_0$ by the
bounded convergence theorem, as required.
\Qed

\procl q.taudecay
The proof actually shows that
$(1/n)\sum_{k=0}^{n-1} \E\Bigl[\tau\bigl(o, X(k)\bigr)\Bigr] \to 0$ as
$n \to\infty$ when there are infinitely many infinite clusters in an
invariant percolation process that has indistinguishable clusters. 
Do we have $\tau\bigl(o,X(n)\bigr)\to 0$ a.s.?
\endprocl

\mnote{Do we have $\E[\tau(o, X_n)] \to 0$?}

\procl x.delete-tau
We give an example of an ergodic invariant
deletion-tolerant percolation process with
infinitely many infinite clusters and with $\tau$ bounded below.
The percolation process $\omega$ will take place on $\Z^3$.
Let $\{a(n) \st n\in\Z\}$ be independent $\{0,1\}$-valued random variables
with $\P[a(n)=1]=1/2$, and let $A$ be the set of
vertices $(x_1,x_2,x_3)\in\Z^3$ with $a(x_1)=1$. 
Let $\eps>0$ be small, and let $z_e$ be i.i.d.\
uniformly distributed in $[0,1]$ indexed by the edges of
$\Z^3$ and independent of the $a(n)$'s.
Let $\omega_1$ be the set of edges $e$ with both
endpoints in $A$ such that $z_e<1-\eps$, and let
$\omega_2$ be the set of edges $e$ with $z_e<\eps$.
Let $\omega_3$ be the set of edges that have no vertex
in common with $\omega_2\cup A$ and satisfy
$z_e<1-\eps$.  Set $\omega:=\omega_1\cup\omega_2\cup\omega_3$.
\par
It is immediate to verify that when $\eps$
is sufficiently small, a.s.\ $\omega$
has a single infinite component whose intersection
with each component of $A$ is infinite, and has
a single infinite component in each component of the
complement of $A$.  The claimed properties of
this example follow easily.
\endprocl

\bsection {Ergodicity of Delayed Random Walk}{s.ergdrw}

The following theorem is not needed for the rest of the paper.
It is presented here because its proof uses some of the ideas from the proof
of \ref l.freq/.

\procl t.ergdel
\procname{Ergodicity of Delayed Random Walk}
\assumptions.
Let $o\in\verts(G)$ be some base point.
Let $(\P,\omega)$ be a $\Gamma$-invariant ergodic bond percolation process on $G$
with infinite clusters a.s.
Also suppose that $\omega$ has indistinguishable infinite clusters a.s.
Let $\Ph_o$ be the joint law of $\omega$ and two-sided delayed simple random walk
on $\omega$, starting at $o$. 
Let $\ev A$ be an event that is $\Gamma$-invariant
and $\SS$-invariant.  Then $\Ph_o[\ev A \mid C(o)\hbox{ is infinite}]$
is either $0$ or $1$.
\endprocl

\proof  
For each $m,n\in\Z$, let $\ev F_m^n$ be the $\sigma$-field of
$\Gamma$-invariant sets generated by
$\omega$ and the random variables $\Seq{W(j)\st j\in[m,n]}$, where
$W(j)$ is the location of the delayed random walk at time $j$.

Let $\ev C$ be the event that $C(o)$ is infinite.
Let $\eps>0$.
  Then there is an $n\in\N$ and an event
$\ev A'\in \ev F_{-n}^n$ such that $\Ph_o[\ev A'\xor \ev A]<\eps \P[\ev C]$,
and 
therefore
$$
\Ph_o[\ev A'\xor\ev A  \mid \ev C] < \eps 
\,.
\label e.apx
$$

Note that $\SS \ev C=\ev C$.  Hence we have for all $m\in\Z$
$$
\Ph_o[\SS^m \A' \xor \ev A  \mid \ev C] < \eps
\,.
$$
Observe also that the two events
$\SS^{n+1} \A'$ and $\SS^{-n-1} \A'$ are independent given $\omega$
by the Markov property for the delayed random walk.
Consequently,
\begineqalno
\int_{\omega\in\ev C} \Ph_o[\ev A \mid\omega]^2 \,d\P[\omega]
&
\geq
\int_{\omega\in\ev C}
\Ph_o[\SS^{n+1} \A'  \mid\omega]
\Ph_o[\SS^{-n-1} \A'  \mid\omega]
\,d\P[\omega]
- 2 \eps
\cr &
=
\int_{\omega\in\ev C} \Ph_o[\SS^{n+1} \A', \SS^{-n-1} \A' \mid\omega]
\,d\P[\omega]
- 2 \eps
\cr &
\geq
\int_{\omega\in\ev C} \Ph_o[\ev A  \mid\omega] \,d\P[\omega]
- 4 \eps
\,.
\endeqalno
Since $\eps$ is arbitrary, it follows that
$$
\int_{\ev C} \Ph_o[\ev A \mid\omega]^2 \,d\P[\omega]
\geq \int_{\ev C} \Ph_o[\ev A \mid\omega] \,d\P[\omega]
\,,
$$
which means that $\Ph_o[\ev A \mid\omega]$ is $0$
or $1$ for almost every $\omega\in\ev C$.

Let $\ev B$ be the set of $\omega\in\ev C$ such that
$\Ph_o[\ev A \mid\omega]=1$ and
let 
$$
\hat{\ev B}:=
\Bigl\{\bigl(C(o),\omega\bigr)\st \omega\in\ev B\Bigr\}
\,.
$$
Finally, let $\Gamma\hat{\ev B}$ denote the orbit of
$\hat{\ev B}$ under $\Gamma$; that is,
$\Gamma\hat{\ev B}:=\{\big(C(\gamma o), \gamma \omega \big) \st \omega\in
{\ev B}\}$.
Because $\SS\ev A=\ev A$ and $\ev A$ is $\Gamma$-invariant,
it follows as in the beginning of the proof of \ref l.freq/ that for every
$\gamma\in\Gamma$ and every $\omega$ such that
$\gamma o\in C(o)$, we have $\omega\in\ev B$ iff
$\gamma^{-1}\omega\in\ev B$ (after possibly making a measure zero
modification of $\ev B$).
Therefore 
$$
\Gamma\hat{\ev B} \cap \{(C,\omega)\st o\in C\} = \hat{\ev B}
\,.
\label e.gamhatb
$$

Since $\omega$ is ergodic and has indistinguishable infinite
clusters, a.s.\ all infinite clusters of $\omega$ are
in $\Gamma\hat{\ev B}$ or a.s.\ all infinite clusters of
$\omega$ are not in $\Gamma\hat{\ev B}$.  If the latter is the
case, then $\P[\ev B]=0$ and hence $\Ph_o[\ev A  \mid \ev C]=0$.
Therefore, assume that
$$
\P[\big(C(o),\omega\big)\in \Gamma \hat{\ev B}  \mid \ev C] =1
\,.
$$
Hence, by \ref e.gamhatb/, $\P[\ev B  \mid \ev C]=1$,
giving $\Ph_o[\ev A  \mid\ev C]=1$, as required.
\Qed

\bsection {Uniqueness for Bernoulli Percolation}{s.uniq}

We now present some applications of \ref t.tau/ to Bernoulli percolation on
Cayley graphs.
Later in this section, we prove inequalities relating $p_u$ based on
\ref t.cerg/.

We first note that the choice of generators does not influence whether $p_u
< 1$:

\procl t.uniqgen
Let $S_1$ and $S_2$ be two finite generating sets for a countable group
$\Gamma$, yielding corresponding Cayley graphs $G_1$ and $G_2$. Then
$p_u(G_1) < 1$ iff $p_u(G_2) < 1$.
\endprocl

\proof
Left and right Cayley graphs with respect to a given set of generators are
isomorphic via $x \mapsto x^{-1}$, so we consider only right Cayley graphs.

Write $\tau_p^i(x)$ for the probability under Bernoulli($p$) percolation that
$o$ and $x$ lie in the same cluster of $G_i$.
Express each element $s \in S_1$ in terms of a word $\varphi(s) \in S_2$.
Let $\omega_2$ be Bernoulli($p$) percolation on $G_2$ and define $\omega_1$
on $G_1$ by letting $[x, xs] \in \omega_1$ iff the path from $x$ to $xs$ in
$G_2$ given by $\varphi(s)$ lies in $\omega_2$.
Then $\omega_1$ is a percolation such that if two edges are sufficiently
far apart, then their presence in $\omega_1$ is independent.
Thus, Liggett, Schonmann and Stacey (1997) provide a function $f(p) \in (0,
1)$ such that $f(p) \uparrow 1$ when $p \uparrow 1$ and such that
$\omega_1$ stochastically dominates Bernoulli$\big(f(p)\big)$ percolation
on $G_1$.
This implies that $\tau_{f(p)}^1(x) \le \tau_p^2(x)$ for each $x$. If $f(p)
{}> p_u(G_1)$, then it follows from this and \ref t.tau/ that $p \ge
p_u(G_2)$, showing that $p_u(G_1) < 1$ implies $p_u(G_2) < 1$.
\Qed

\procl r.nolss
For the situation used in the above proof, and for many similar applications,
one does not need the full generality of the theorem of Liggett,
Schonmann and  Stacey.  The following observation suffices.
Suppose that $(X_i\st i\in I)$ are i.i.d.\ random variables
taking values in $\{0,1\}$ and with $\P[X_i = 1]=p$. Let
$(I_j\st j\in J)$ be an indexed collection of subsets of $I$.
Set $Y_j:=\min\{X_i\st i\in I_j\}$ and $J_i:=\{j\in J\st i\in I_j\}$.
Suppose that $n:=\sup\{|I_j|\st j\in J\}$ and $m:=\sup\{|J_i|\st i\in I\}$
are both finite.
Then $(Y_j\st j\in J)$ stochastically dominates independent
random variables $(Z_j\st j\in J)$ with
$\P[Z_j=1]\geq \left(1-(1-p)^{1/m}\right)^n$.
Indeed, set $X'_i:= \max\{ Z_{i,j}\st j\in J_i\}$
and $Z_j := \min\{Z_{i,j}\st i\in I_j\}$, where 
$(Z_{i,j}\st i\in I,\ j\in J)$ are independent $\{0,1\}$-variables
with $\P[Z_{i,j}=1] =1-(1-p)^{1/m}$.
Then $(X_i\st i\in I)$ stochastically dominates
$(X'_i\st i\in I)$, and hence $(Y_j\st j\in J)$
stochastically dominates $(Z_j\st j\in J)$.
\endprocl

\procl r.puinvgen
Schonmann~(1999a) shows that for every quasi-transitive graph $G$,
$p_u=p_{BB}$, where $p_{BB}$ is the infimum of all $p$ such that 
$$
\lim_{r\to\infty}\inf_{v,u\in\verts}
\P\bigl[B(v,r)\leftrightarrow B(u,r)\bigr]=1
$$
in Bernoulli($p$) percolation, where $B(v,r)$ denotes the ball in $G$ of
radius $r$ and center $v$ and $A\leftrightarrow A'$ is the event
that there is a cluster $C$ with $C\cap A\neq\emptyset$ and
$C\cap A'\ne\emptyset$.  Based on this and the proof
of \ref t.uniqgen/, one obtains the following generalization.
Suppose that $G$ and $G'$ are quasi-transitive graphs and 
$G'$ is quasi-isometric to $G$. Then $p_u(G)<1$ iff $p_u(G')<1$.
This observation was also made independently by Y.~Peres (private
communication).
\endprocl

For our next result regarding $p_u$, we need the following construction.
Let $\kappa_0$ be a real-valued random variable, and suppose that
$\P$ is a bond percolation process on some graph $G$.
We would like to color the clusters of $\P$ in such a way
that conditioned on the configuration $\omega$, the colors of the components
are i.i.d.\ random variables with the same law as $\kappa_0$.

To construct this process, let $\Seq{v_1,v_2,\dots}$
be an ordering of the vertices in $G$.  Let $\Seq{\kappa_1,\kappa_2,\dots}$
be i.i.d.\ random variables with the same law as $\kappa_0$.
Given $\omega$ and $v\in\verts$, set $\kappa_\omega(v):=\kappa_j$
if $j$ is the least integer with $v_j\in C(v)$.
It is not hard to see that $\kappa_\omega(v)$ is measurable, as a function
on the product of $2^\edges$ and the sample spaces of the $\kappa_j$'s.
Let $\P^\kappa$ denote the law of $(\omega,\kappa_\omega)$.
Observe that $\P^\kappa$ satisfies the description of the
previous paragraph, and therefore does not depend on the choice
of the ordering of $\verts$. 
Consequently, if $\gamma$ is an automorphism of $G$
and $\P$ is $\gamma$-invariant, then $\P^\kappa$ is $\gamma$-invariant.

\procl l.colored
\procname{Ergodicity of $\P^\kappa$}
\genassumptions.
Suppose that $(\P,\omega)$ is a $\Gamma$-invariant ergodic insertion-tolerant
percolation process on $G$.
Let $\kappa_0$ be a real-valued random variable, and let $\P^\kappa$
be as above.
Suppose that $\inf\bigl\{ \tau(o,x)\st {x \in \verts(G)}\bigr\}= 0$.
Then $\P^\kappa$ is $\Gamma$-invariant and ergodic.
\endprocl

\proof
To prove the ergodicity of $\P^\kappa$,
let $\ev A$ be a $\Gamma$-invariant event in $2^{\edges(G)} \times \R^{\verts(G)}$
and $\epsilon \in(0,1/2)$. 
The probability of $\A$ conditioned on $\omega$ must
be a constant, say $a$, by ergodicity of $\P$.  We need to show
that $a \in \{0, 1\}$.

There is a cylindrical event $\ev A'$ 
that depends only on the restriction of
$\omega$ and $\kappa_\omega$ to some ball $B(o,r)$ about $o$ and
such that $\P^\kappa[\ev A \xor \ev A'] < \epsilon$. 
Then $\EBig{\big|\P^\kappa[\A' \mid \omega] - a\big|} < \epsilon$.

Let $\ev B_x$ be the event that some vertex in $B(o,r)$ belongs to the same
cluster as some vertex in $B(x,r)$.
Since $\inf_x \tau(o,x) = 0$ and $\P$ is insertion tolerant,
there is some $x$ such that $\P[\ev B_x] < \epsilon$.
Indeed, suppose that $\P[\ev B_x] \ge \epsilon$ for all $x$.
Let $\ev D_x$ be the event that all the edges in $B(x, r)$ belong
to $\omega$, and for $A\subset G$, let $\F_{\neg A}$ denote the
$\sigma$-field generated by the events $\{e\in\omega\}$ with $e\notin A$.
Then $\ev B_x$ is $\F_{\neg (B(o,r)\cup B(x,r))}$-measurable.
By insertion tolerance, for every $\F_{\neg B(o,r)}$-measurable
event $\ev C$ with $\P[\ev C]>0$, we have $\P[\ev C\cap \D_o]>0$. 
Consequently, there is some $\delta>0$ such that
$$
\P\Bigl[\P\bigl[\D_o\bigm|\F_{\neg B(o,r)}\bigr]>\delta\Bigr]>1-\eps/2
\,.
$$
It follows that for all $\F_{\neg B(o,r)}$-measurable events $\ev C$ with
$\P[\ev C]\geq\eps$, we have
$\P\bigl[\D_o\bigm|\ev C\bigr]>\delta/2$.
In particular, for all $x\in\verts$, we have
$\P\bigl[\D_o\bigm|\ev B_x\bigr]>\delta/2$,
which gives
$\P\bigl[\D_o\cap\ev B_x\bigr]>\eps\delta/2$.
Transitivity and insertion tolerance imply that there is a $\delta'>0$
such that
$$
\P\Bigl[\P\bigl[\D_x\bigm|\F_{\neg B(x,r)}\bigr]>\delta'\Bigr]>1-\eps\delta/4
$$
for all $x\in\verts$.
Hence, for all $\F_{\neg B(x,r)}$-measurable events $\ev C$ with
$\P[\ev C]\geq \eps\delta/2$, we have
$\P\bigl[\D_x\bigm|\ev C\bigr]>\delta'/2$.
Taking $\ev C:= \D_o\cap\ev B_x$ gives for all $x\in\verts\setminu B(o,2r)$
$$
\tau(o, x) \ge \P[\ev \D_o \cap \ev B_x \cap \ev \D_x]
=\P[\ev \D_x\mid \ev \D_o\cap\ev B_x]
\P[\ev \D_o\cap\ev B_x]
\ge (\delta'/2)(\eps\delta/2)
\,,
$$
which contradicts our assumption, and thereby verifies that
$\P[\ev B_x] < \epsilon$ for some $x\in \verts$.  Fix such an $x$.

Let $\gamma_x\in\Gamma$ be such that $\gamma_x o = x$.
Then 
$$
\EBig{\big|\P^\kappa[\gamma_x \A' \mid \omega] - a\big|} < \epsilon
\,.
$$
Since $\gamma_x \A'$ depends only on the colored configuration in $B(x, r)$, we
have that on the complement of $\ev B_x$,
$$
\P^\kappa[\A', \gamma_x \A'\mid \omega] = \P^\kappa[\A' \mid \omega] \P^\kappa[\gamma_x
\A' \mid \omega]
\,,
$$
whence
$\EBig{\big|\P^\kappa[\A', \gamma_x \A'\mid \omega] - a^2\big|} = O(\epsilon)$.
Therefore, $a - a^2 =
\EBig{\big|\P^\kappa[\A, \gamma_x \A\mid \omega] - a^2\big|} = O(\epsilon)$.
Since $\epsilon$ was arbitrary, we get that $a \in \{0, 1\}$, as desired.
\Qed

We now prove that $p_u < 1$ for all Cayley graphs of
Kazhdan groups.
Let $\Gamma$ be a countable group and $S$ a finite subset of $\Gamma$.  Let
$\U(\H)$ denote the set of unitary representations of $\Gamma$ on a Hilbert
space $\H$ that have no invariant vectors except $0$. Set
$$
\kappa(\Gamma, S) := \max \Bigl\{ \epsilon \st \forall \H \ \forall \pi \in
\U(\H)\ 
\forall v \in \H \ \exists s \in S \quad \|\pi(s)v - v\| \ge \epsilon \|v\|
\Bigr\}\,.
$$
Then $\Gamma$ is called {\bf Kazhdan} (or has Kazhdan's property (T)) if
$\kappa(\Gamma, S) > 0$ for all finite $S$. 
The only amenable Kazhdan groups
are the finite ones. Examples of Kazhdan groups include SL$(n, \Z)$ for $n
\ge 3$.  See de la Harpe and Valette (1989) for background;
in particular, every Kazhdan group is finitely generated (p.~11), but not
necessarily finitely presentable (as shown by examples of Gromov; see p.~43).
It can be shown directly, but also follows from our \ref c.kazhdan/ below,
that every infinite Kazhdan group has only one end.
See Zuk (1996) for examples of Kazhdan groups arising as fundamental
groups of finite simplicial complexes.

Rather than the definition, we will use the following characterization of
Kazhdan groups:  Let $\P_*$ be the
probability measure on subsets of $\Gamma$ that is the empty set half the time
and all of $\Gamma$ half the time.  Recall that $\Gamma$ acts by
translation on the probability measures on $2^\Gamma$.

\procl t.GW \procname{Glasner and Weiss (1998)}
A countable infinite group $\Gamma$ is Kazhdan iff $\P_*$ is not in the
weak${}^*$ closure of the $\Gamma$-invariant ergodic probability measures
on $2^\Gamma$.
\endprocl

\procl c.kazhdan
If $G$ is a Cayley graph of an infinite Kazhdan group $\Gamma$, then
$p_u(G) < 1$. Moreover, $\P[\exists\hbox{ a unique infinite cluster}]=0$
in Bernoulli($p_u$) percolation.
\endprocl

In an older version of this manuscript, only the first statement appeared.
We thank Yuval Peres for pointing out that a modification
of the original proof produces the stronger second statement.

Schonmann (1999b) proved that Bernoulli($p_u$) on $T\times \Z$, where
$T$ is a regular tree of degree at least three, does not have a unique
infinite cluster, and Peres (1999) generalized this result to 
nonamenable products.  On the other hand, Bernoulli($p_u$) percolation
on a planar nonamenable transitive graph has a unique infinite cluster.
(See Lalley (1998) for the high genus case, and Benjamini and Schramm
(in preparation) for the general case.)  Little else is known, however,
about the uniqueness of infinite clusters at $p_u$. For
example, the case of lattices in
hyperbolic $3$-space is still open.

\proof
Suppose that $p_u(G) = 1$.
Let $o$ be the identity in $\Gamma$, regarded as a vertex in $G$.
Write $\tau_p(x)$ for the probability under Bernoulli($p$) percolation that
$o$ and $x$ lie in the same cluster.
Then by \ref t.tau/, for all $p<1$ we have $\inf_x \tau_p(x) = 0$.
Fix $p$ and let $\omega$ be the open subgraph of a Bernoulli($p$)
percolation. 
Let $\eta$ be the union of the sites of some of the clusters of $\omega$,
where each cluster is independently put in $\eta$ with probability 1/2. 
By \ref l.colored/, the law $\law$ of $\eta$ is $\Gamma$-invariant and ergodic.
Furthermore, any fixed finite subset of $\Gamma$ either is contained in
$\eta$ or is disjoint from $\eta$ with high probability when $p$ is
sufficiently close to 1.
That is, $\P_*=\hbox{weak}^*\hbox{-}\lim_{p\to 1}\law$,
whence $\Gamma$ is not Kazhdan.

To prove the stronger statement, suppose that there is a unique infinite
component $\P$-a.s. Let $(\omega, \widetilde\omega)$ be the standard
coupling of Bernoulli($p)$ and Bernoulli($p_u$) percolation on $G$.
Let $\eta$ be as above. For a vertex $x \in G$, write $A(x)$ for the set of
clusters of $\omega$ that lie in the unique infinite cluster of
$\widetilde\omega$ and that are closest to $x$ among those with this
property (where distance is measured in $G$). Define $\eta'$ to be the
union of $\eta$ with all sites $x$ for which $A(x)$ contains only one
cluster and that cluster lies in $\eta$.
Since the law of $(\omega, \widetilde\omega, \eta)$ is ergodic by an
obvious extension of \ref l.colored/, so is the
law $\law'$ of its factor $\eta'$. Again, we obtain
$\P_*=\hbox{weak}^*\hbox{-}\lim_{p\to p_u}\law'$, whence $\Gamma$ is not
Kazhdan.
\Qed

\procl r.alt-proof
For probabilists, we believe that the proof we have presented of \ref
c.kazhdan/ is
the most natural. For others, we note that one can avoid \ref l.colored/
and \ref t.GW/ by using a theorem of 
Delorme (1977) and Guichardet (1977) that
characterizes Kazhdan groups in terms of positive semidefinite functions.
This relies on the fact that $\tau(x, y)$ is positive semidefinite, as
observed by Aizenman and Newman (1984).
\endprocl

Other groups that are not finitely presentable and that have provided
interesting examples for probability theory are the so-called lamplighter
groups (see Kaimanovich and Vershik (1983) and Lyons, Pemantle and Peres
(1996); in these references, these groups are denoted $G_d$ and are
amenable, but we will be interested here in nonamenable examples). 

We first give a concrete description of a {\bf lamplighter graph},
and later generalize and use more algebraic language.
Suppose that $G$ is a graph.  The {\bf lamplighter graph} $L_G$
over $G$ is the graph whose vertices are pairs $(A,v)$,
where $A\subset\verts(G)$ is finite and $v\in\verts(G)$.
We think of $A$ as the locations of the lamps that are on, and
consider $v$ as the location of the lamplighter.
One neighbor of $(A,v)$ in $L_G$ is the vertex $(A\xor\{v\},v)$
(the lamplighter switches the lamp off or on)
and the other neighbors have the form $(A,u)$, where $[v,u]\in\edges(G)$
(the lamplighter walks one step).

In the algebraic context and language,
lamplighter groups are particular wreath products: Let $\Gamma$ be a
group acting from the left on a set $\verts$.
Let $K$ be a group; $K^\verts_*$ denotes the group of maps $f : \verts \to
K$ such that $f(x)$ is the identity element $\id_K$ of $K$ for all but finitely
many $x \in \verts$ and with multiplication $(f_1 f_2)(x) := f_1(x)
f_2(x)$.
Then $\Gamma$ acts on $K^\verts_*$ by translation: $(\gamma
f)(x) := f(\gamma^{-1} x)$.
The {\bf (restricted) wreath product} $K \wreath \Gamma$
is the set $K^\verts_* \times \Gamma$ with the multiplication
$$
(f_1, \gamma_1) (f_2, \gamma_2) 
:=
\big(f_1 (\gamma_1 f_2), \gamma_1 \gamma_2\big)
\,.
$$

If $\Gamma$ and $K$ are finitely generated and $\Gamma$ acts transitively on $V$,
then $K \wreath \Gamma$ is finitely generated. 
To see this, let $\gamma_1, \ldots, \gamma_s$ generate $\Gamma$ and $k_1,
\ldots, k_t$ generate $K$. 
Write $\id_\Gamma$ for the identity element of $\Gamma$ and $\id_\verts$
for the identity element of $K^\verts_*$.
Let $o \in \verts$ be fixed.
Write $F_j$ for the element of $K^\verts_*$ defined by $F_j(o) := k_j$ and
$F_j(x) := \id_K$ for all $x \ne o$.
Set 
\begineqalno
S_\Gamma:= & \bigl\{(\id_\verts, \gamma_i)\st 1\le i\le s\bigr\}\,,\cr
S_K:= & \bigl\{ (F_j, \id_\Gamma)\st 1 \le j \le t \bigr\}
\,.
\endeqalno
Then $S:=S_\Gamma\cup S_K$ is a finite generating set for $K \wreath \Gamma$.
Indeed, let $(f, \gamma) \in K \wreath \Gamma$.
For $x \in \verts$, choose $\gamma_x \in \Gamma$ such that
$\gamma_x o = x$ and write $h_x \in K^\verts_*$ for the function $h_x(o) :=
f(x)$ and $h_x(y) := \id_K$ for $y \ne o$. Then
$$
(f, \gamma) =
\left(\prod_{x \in \verts \atop f(x) \ne \id_K}
(\id_\verts, \gamma_x) (h_x, \id_\Gamma) (\id_\verts, \gamma_x^{-1})
\right)
(\id_\verts, \gamma)
\,,
\label e.path
$$ %sz
where the product over $x$ is taken in any order.
If we then write each $(h_x, \id_\Gamma)$ as a
product of $(F_j, \id_\Gamma)$, each $(\id_\verts, \gamma_x)$ as a product
of $(\id_\verts, \gamma_i)$, each $(\id_\verts, \gamma_x^{-1})$ as a product
of $(\id_\verts, \gamma_i)$, and $(\id_\verts, \gamma)$ as a product
of $(\id_\verts, \gamma_i)$, we obtain a representation as a product of elements
of $S$.

The lamplighter groups $G_d$ are those where $\Gamma = \verts = \Z^d$
and $K = \Z_2$.
(Rick Kenyon (private communication) has observed that the groups $G_d$ are not
finitely presentable.)

\procl c.wreath
Let $K$ be a finite group
with more than one element.
Let $\Gamma$ be a finitely
generated group acting transitively on an infinite set $\verts$. 
If $G$ is any Cayley graph of $K \wreath \Gamma$,
then $p_u(G) < 1$.
\endprocl

To prove this, we borrow a technique from Benjamini, Pemantle and
Peres (1998):
%If $\psi$ is a path, denote by $|\psi|$ its length. 
If $\phi$ and $\psi$ are two paths, denote by $|\phi \cap \psi|$ the
number of edges they have in common (as sets of edges).

\procl l.EIT
Let $G$ be a graph and $x, y$ be any two vertices in $G$.
Let $\theta \in (0, 1)$ and $c > 0$ be constants.
Suppose that $\mu$ is a probability measure on
(possibly self-intersecting) paths $\psi$ joining $x$ to $y$ such that 
$$
(\mu\times\mu)\{ (\phi, \psi) \st |\phi \cap \psi| = n\} \le c \theta^n
$$
for all $n \in \N$.
Let $\theta < p < 1$.
Then for Bernoulli($p$) percolation on $G$, we have
$$
\tau(x, y) \ge c^{-1} (1 - \theta/p)
\,.
$$
\endprocl

\proof
Define the random variable 
$$
Z := \sum_{\psi} \mu(\psi) {\II{\psi\rm\ is\ open} \over \P[\psi \hbox{ is
open}]}
\,.
$$
Then $\E[Z] = 1$ and
\begineqalno
\E[Z^2]
&=
\sum_{\phi,\psi} \mu(\phi) \mu(\psi) {\P[\phi, \psi \hbox{ are both open}]
\over \P[\phi \hbox{ is open}] \P[\psi \hbox{ is open}]}
\cr&=
\sum_{\phi,\psi} \mu(\phi) \mu(\psi) p^{-|\phi \cap \psi|} 
\cr&=
\sum_n p^{-n} \sum_{|\phi \cap \psi| = n} \mu(\phi) \mu(\psi) 
\cr&\le
\sum_n p^{-n} c \theta^n
=
c (1 - \theta/p)^{-1}
\,.\cr
\endeqalno
Therefore, the Cauchy-Schwarz inequality yields
$$
\tau(x, y) \geq \P[Z > 0] \ge \E[Z]^2 / \E[Z^2] \ge c^{-1} (1 - \theta/p)
\,. \Qed
$$

\proofof c.wreath
Since $K$ is assumed to be finite, we take the generating set
$\{k_1,\dots,k_t\}$ to be all of $K$.
We use the Cayley graph given by the generating set $S$.
It suffices to exhibit a measure on paths connecting $(\id_{\verts},
\id_\Gamma)$ to $(f, \gamma)$ that satisfies the condition of \ref l.EIT/
with $\theta$ and $c$ not depending on $(f, \gamma)$, since then \ref t.tau/
implies that $p_u \le \theta$.

Define edges joining pairs in $\verts$ by
$$
\edges := \{[x, \gamma_i^\epsilon x] \st x \in \verts,\; 1 \le i \le s,\;
\epsilon = \pm 1 \}
\,.
$$%stopzone
The graph $G' := (\verts,\edges)$ will be called the {\bf base graph}.
Note that the base graph is not the same as the Cayley graph,
$G=(\verts(G),\edges(G))$.
Let $\pi:\verts(G)\to\verts$ be the projection $\pi(g,\beta):=\beta o$.
Let $v_1,v_2,\dots$ and $u_1,u_2,\dots$ be 
infinite simple paths in $G'$ starting at $v_1=o$ and $u_1=\gamma o$,
respectively, that are disjoint, except that $v_1$ may equal $u_1$.
Because $G'$ is an infinite, connected, transitive graph,
it is easy to show that such paths exist.
For each $j=1,2,\dots$,
fix an $\alpha_j\in S_\Gamma$ such that $\alpha_j v_j=v_{j+1}$
and fix a $\beta_j\in S_\Gamma$ such that $\beta_j u_j=u_{j+1}$.

Given a word $W=w_1w_2\cdots w_n$ with letters from $S$, let
$W^{-1}$ denote the word $w_n^{-1}w_{n-1}^{-1}\cdots w_1^{-1}$,
and let $\tilde W(j)$, $j\in\{0,1,\dots,n\}$, denote the group
element $w_1w_2\cdots w_j$.

Let $W_\gamma$ be a word in $S_\Gamma$ representing $(\id_\verts,\gamma)$
with the property that for every $v\in\verts$ such that
$f(v)\neq\id_K$, there is a $j$
such that $\pi\tilde W(j)=v$.  Let $n$ be the length of $W_\gamma$.
Let $W_\alpha$ be the word $\alpha_1\alpha_2\cdots\alpha_n$,
and let $W_\beta$ be the word $\beta_1\beta_2\cdots\beta_n$.
Let $W$ be the concatenation
$$
W:=W_\alpha W_\alpha^{-1}W_\gamma W_\beta W_\beta^{-1} W_\gamma ^{-1} 
W_\alpha W_\alpha^{-1}W_\gamma W_\beta W_\beta^{-1} \,,
% \label e.con
$$
and let $N$ be the length of $W$.
Let the letters in $W$ be $W=w_1w_2\cdots w_N$.
Consider words of the form $\phi(X):= w_1 X_1 w_2 X_2\cdots w_N X_N$,
where $X=(X_1,\dots,X_N)\in S_K^N=(S_K)^N$.
For any $v\in\verts$, let $J_v$ be the set of $j\in\{1,2,\dots,N\}$
such that $\pi\tilde W(j) =v$. 
Let $\pi_K : K \wreath \Gamma \to K^\verts_*$ be the projection onto the first
coordinate.
The uniform measure $\mu_0$ on the set of $X\in S_K^N$
such that $\phi(X)$ is a word
representing $(f,\gamma)$ can be described as follows: 
$\prod_{j\in J_v} \pi_K X_j(o)=f(v)$
(where the order of multiplication is the order of $J_v$ as a subset of $\N$)
and for every $j_1\in J_v$, the random variables
$\bigl(X_j\st j\in J_v,\ j\neq j_1\bigr)$ are independent,
uniform in $S_K$, and independent of $(X_j\st j\notin J_v)$.
Note that $\phi(X)$ can also be thought of as a random path 
in $G$ from $(\id_{\verts}, \id_\Gamma)$ to $(f, \gamma)$.

Let $X$ and $Y$ be i.i.d.\ with law $\mu_0$.
We want to bound the probability that
there are $k$ edges shared by the paths $\phi(X)$ and $\phi(Y)$.
In fact, we shall bound the
probability that these paths share at least $k$ vertices.
% Consider the description \ref e.con/ of $W$. 
Given any $j\in\{1,\dots,N\}$, let $H(j)$ be the
set of $v\in \verts$ such that
$\min J_v < j < \max J_v$,
and let $h(j):=\left| H(j)\right|$.
The choice of $W$ ensures that $h(j)\geq \min\{j,N-j,n\}-1$.
Because $N=O(n)$, this gives
$$
h(j)\geq c_1\min\{j,N-j\}-1
\,,
$$
for some universal constant $c_1>0$.
If $(g,\delta)$ is the element of $\hat\verts$
represented by the word $w_1 X_1 w_2\cdots w_j X_j$,
then $\bigl(g(v)\st v\in H(j)\bigr)$, are i.i.d.\ uniform
in $K$.  Consequently, the probability that
$w_1 X_1 w_2\cdots w_j X_j =  w_1 Y_1 w_2\cdots w_{j'} Y_{j'}$,
as elements of $K \wreath \Gamma$, is at most
$|K|^{-\max\{h(j),h(j')\}}$.
If there are more than $8k$ vertices common 
to $\phi(X)$ and $\phi(Y)$, then there must be a pair of indices
$j,j'\in\{k,k+1,\dots,N-k\}$ such that
$w_1 X_1 w_2\cdots w_j X_j =  w_1 Y_1 w_2\cdots w_{j'} Y_{j'}$,
as elements of $K \wreath \Gamma$, or with a similar 
equality when the rightmost letter is dropped from either
or both sides.
The probability for that is at most
$$
4\sum_{j=k}^{N-k}\sum_{j'=k}^{N-k}
|K|^{-\max\{h(j),h(j')\}}
\leq
4\sum_{j=k}^{N-k}\sum_{j'=k}^{N-k}
|K|^{-(h(j)+h(j'))/2}
\leq c_2 |K|^{-c_1 k}
\,,
$$
% by our estimate for $h(j)$,
where $c_2$ is a constant depending only on $|K|$. 
Consequently, an appeal to \ref l.EIT/ completes the proof.
\Qed

Benjamini and Schramm (1996) conjectured that for any nonamenable Cayley
graph $G$, we have $p_c(G) < p_u(G)$. This is still open, but using our
main result on cluster indistinguishability, we may show that a comparable
statement fails for invariant percolation processes:

\procl c.phases
There is an invariant deletion-tolerant percolation process
$(\P,\omega)$ on a nonamenable Cayley
graph that has only finite components a.s., but such that
for every $\epsilon>0$, the union of $\omega$ with an independent 
Bernoulli($\epsilon$) percolation process produces a unique infinite cluster
a.s.
\endprocl

\proof 
Recall that Bernoulli($p$) bond percolation on $\Z^2$
produces a.s.\ no infinite cluster iff $p\leq 1/2$ (see Grimmett 1989).
Let ${\Bbb F}_2$ be the Cayley graph of the free group on two letters with
the usual generating set.
Let $G:={\Bbb F}_2\times \Z^2$.  Each edge of $G$ joins vertices 
$(x_1, y_1)$ to $(x_2, y_2)$, where either $x_1 = x_2$ and $y_1 \sim y_2$,
or $x_1 \sim x_2$ and $y_1 = y_2$.
The former type of edge will be called a $\Z^2$ edge, while the latter
will be called an ${\Bbb F}_2$ edge.
For each $x \in {\Bbb F_2}$, the set of all edges $[(x, y_1), (x, y_2)]$ is
called the $\Z^2$ fiber of $x$.
Let $\omega$ be Bernoulli($1/2$) bond percolation on all the $\Z^2$
edges; it contains no ${\Bbb F}_2$ edge.
Then $\omega$ is invariant and deletion tolerant.
When we take the union with an independent Bernoulli($\epsilon$) process
$\eta$, the intersection of $\omega \cup \eta$ with each $\Z^2$ fiber
contains a.s.\ a unique infinite cluster. 
These clusters a.s.\ connect to each other in $G$. 
Thus, $\omega \cup \eta$ contains a unique infinite cluster whose
intersection with each $\Z^2$ fiber has an infinite component.
Since $\omega \cup \eta$ is insertion tolerant, 
cluster indistinguishability of $\omega \cup \eta$ implies that it has no
other infinite cluster.
\Qed

\procl q.iuniq
Is there an insertion-tolerant process with the property exhibited in \ref
c.phases/?
\endprocl

The same line of thought is also used to prove the following theorem,
which is an answer to a question posed by Yuval Peres (private communication)
and motivated by the work of Salzano and Schonmann (1997, 1998)
on contact processes.

\procl t.puineq 
Let $G$, $H$ and $H'$ be unimodular transitive graphs.
Assume that $H'\subset H$ and $G$ is infinite.  Then
$$
p_u(G\times H')\geq p_u(G\times H) 
\,.
$$
In particular, when $H'$ consists of a single vertex of $H$ we get
$$
p_u(G) \geq p_u(G\times H)
\,.
$$
\endprocl %sz

\proof
Let $\omega$ be Bernoulli($p$) bond percolation on $G\times H$,
where $p>p_u(G\times H')$, and let $\Gamma$ be the product
of a unimodular transitive automorphism group on $G$
with a unimodular transitive automorphism group on $H$.
Then $\Gamma$ is a unimodular transitive automorphism group on
$G\times H$.

Let ${\hat Z}:=\{\gamma(G\times H')\st \gamma\in\Gamma\}$.
By \ref t.uniq/ applied to $G\times H'$, 
for every $Z\in{\hat Z}$ a.s.\ there is a unique
infinite cluster, say $Q_Z$, of $\omega\cap Z$.
Now let $Z,Z'\in{\hat Z}$.  We claim that a.s.\
$Q_Z$ and $Q_{Z'}$ are in the same cluster of $\omega$.
Indeed, since $G$ is infinite, there are
infinite sequences of vertices $v_j\in Z$, $v_j'\in Z'$
such that $\dist(v_j,v_j')=\inf\{\dist(v,v')\st v\in Z,\ v'\in Z'\}$
for all $j$. 
Because $\inf_j\P[v_j\in Q_Z ,\ v_j'\in Q_{Z'}]>0$,
with positive probability
$v_j\in Q_Z$ and $v_j'\in Q_{Z'}$ for infinitely many $j$,
and by Kolmogorov's 0-1 law this holds a.s. 
It is then clear that for some such $j$ there is a connection
in $\omega$ between $v_j$ and $v_j'$.
It follows that a.s.\ there is a unique cluster $Q$ of
$\omega$ that contains every $Q_Z$, $Z\in{\hat Z}$.

Let $\ev A$ be the set of all pairs $(C,W)\in 2^{\verts(G\times H)}
\times 2^{\edges(G\times H)}$ such that
there is a $Z\in\hat Z$ so that $C$ meets an infinite component
of $Z\cap W$.  We know that $Q$ is the only cluster of
$\omega$ with $(Q,\omega)\in\ev A$.  By \ref t.cerg/, 
it follows then that $Q$ is the only infinite cluster of
$\omega$.
\Qed

\procl r.genpuineq 
The same theorem holds without the assumption of unimodularity
and with ``transitive'' replaced by ``quasi-transitive''.
This follows similarly from the indistinguishability of robust properties
proved by H\"aggstr\"om, Peres and Schonmann (1998).
\endprocl

\procl r.nogen
There are infinite Cayley graphs $H'\subset H$ with
$p_u(H')< p_u(H)$.  For example, one may take $H'=\Z^2$
and let $H$ be the free product $\Z^2*\Z_2$, as in \ref r.directions/.
\endprocl 

The work of Schonmann (1999a) was motivated by an analogy between
percolation and contact processes.
More specifically, the property of having complete convergence
with survival for a contact process is closely analogous to
having uniqueness of the infinite cluster for percolation.
However, \ref r.nogen/ describes an instance where this analogy fails,
because complete convergence with survival is a property which is
monotone in the graph (Salzano and Schonmann 1997).

\medskip\noindent{\bf Acknowledgement.} We are grateful to Itai Benjamini,
Olle H\"aggstr\"om, Yuval Peres, and Roberto Schonmann for fruitful
conversations and helpful advice.

\beginreferences

Aizenman, M., Kesten, H. \and Newman, C.M. (1987) Uniqueness of the
infinite cluster and continuity of connectivity functions for short- and
long-range percolation, {\it Commun. Math. Phys.} {\bf 111}, 505--532.

Aizenman, M. \and Newman, C.M. (1984) Tree graph inequalities and critical
behavior in percolation models, {\it J. Stat. Phys.} {\bf 36}, 107--143.

Alexander, K. (1995) Simultaneous uniqueness of infinite clusters in
stationary random labeled graphs, {\it Commun. Math. Phys.} {\bf 168},
39--55.

Babson, E. \and Benjamini, I. (1999) Cut sets and normed cohomology with
applications to percolation, 
{\it Proc. Amer. Math. Soc.}, 589--597.

Benjamini, I., Lyons, R., Peres, Y., \and Schramm, O. (1999)
Group-invariant percolation on graphs, {\it Geom. Funct. Anal.}, to appear.

Benjamini, I., Lyons, R., \and Schramm, O. (1998) Percolation perturbations
in potential theory and random walks, in {\it Random Walks and Discrete
Potential Theory (Cortona, 1997), Sympos. Math., XXXVII (?)}, Cambridge
Univ. Press, Cambridge, 1998, to appear.

Benjamini, I., Pemantle, R. \and Peres, Y. (1998)
Unpredictable paths and percolation, {\it Ann. Probab.} {\bf 26},
1198--1211.

Benjamini, I. \and Schramm, O. (1996) Percolation beyond $\Z^d$, many
questions and a few answers, {\it Electronic Commun. Probab.} {\bf 1},
71--82.

Burton, R.M. \and Keane, M. (1989)
Density and uniqueness in percolation, {\it Commun. Math.
Phys.} {\bf 121}, % no. 3,
501--505.
% 90g:60090

Delorme, P. (1977) 1-cohomologie des repr\'esentations unitaires des groupes
de Lie semi-simples et r\'esolubles. Produits tensoriels continus et
repr\'esentations, {\it Bull. Soc. Math. France} {\bf 105}, 281--336.

Gandolfi, A., Keane. M.S. \and Newman, C.M. (1992) Uniqueness of the
infinite component in a random graph with applications to percolation and
spin glasses, {\it Probab. Theory Relat. Fields} {\bf 92}, 511--527.

Glasner, E. \and Weiss, B. (1998) Kazhdan's property T and the
geometry of the collection of invariant measures, {\it Geom. Funct. Anal.},
to appear.

Grimmett, G.R. (1989) {\it Percolation}. Springer, New York.

Grimmett, G.R. \and Newman, C.M. (1990) Percolation in $\infty + 1$
dimensions, in {\it Disorder in Physical Systems}, G.R. Grimmett and D.J.A.
Welsh (editors), pp.~219--240. Clarendon Press, Oxford.

Guichardet, A. (1977) \'Etude de la 1-cohomologie et de la topologie du
dual pour les groupes de Lie \`a radical ab\'elien, {\it Math. Ann.} {\bf
228}, 215--232.

H\"aggstr\"om, O. (1997) Infinite clusters in dependent automorphism
invariant percolation on trees, {\it Ann. Probab.} {\bf 25}, 1423--1436.

H\"aggstr\"om, O. \and Peres, Y. (1999) Monotonicity of uniqueness for
percolation on Cayley graphs: all infinite clusters are born
simultaneously, 
{\it Probab. Theory Rel. Fields} {\bf 113}, 273--285.

H\"aggstr\"om, O., Peres, Y., \and Schonmann, R.H. (1998) 
Percolation on transitive graphs as a coalescent
process: relentless merging followed by simultaneous uniqueness,
in {\it Perplexing Probability Problems: Papers in Honor of H.~Kesten},
M. Bramson and R. Durrett (editors). Birkh\"auser, Boston, to appear.

de la Harpe, P. \and Valette, A. (1989) La propri\'et\'e $(T)$ de Kazhdan
pour les groupes localement compacts (avec un appendice de Marc Burger),
{\it Ast\'erisque} {\bf 175}. 

Kaimanovich, V.A. \and Vershik, A.M. (1983) Random walks on
discrete groups: boundary and entropy, {\it Ann. Probab.} {\bf 11},
457--490.

Lalley, S.P. (1998) Percolation on Fuchsian groups,
{\it Ann. Inst. H. Poincar\'e Probab. Statist.} {\bf 34}, 151--177.

Liggett, T.M., Schonmann, R.H. \and Stacey, A.M. (1997) Domination by
product measures, {\it Ann.  Probab.} {\bf 25}, 71--95.

Lyons, R. (1990) Random walks and percolation on trees, {\it Ann. Probab.}
{\bf 18}, 931--958.

Lyons, R., Pemantle, R. \and Peres, Y. (1996) Random walks on
the Lamplighter Group, {\it Ann. Probab.} {\bf 24}, 1993--2006. 

Lyons, R. \and Peres, Y. (1998) {\it Probability on Trees and Networks}.
Cambridge University Press, in preparation. Current version available 
at \hfill\break \htmlref{http://php.indiana.edu/\string~rdlyons/}{\tt http://php.indiana.edu/\string~rdlyons/}.

Lyons, R. \and Schramm, O. (1998) Stationary measures for random
walks in a random environment with random scenery, {\it in preparation}.
%{\tt http://php.indiana.edu/\~{}rdlyons/}.

Newman, C.M. \and Schulman, L.S. (1981) Infinite clusters in percolation
models, {\it J. Statist. Phys.} {\bf 26}, 613--628.

Peres, Y. (1999) Percolation on nonamenable products at the uniqueness
threshold, {\it preprint}.

Salzano, M. \and Schonmann, R.H. (1997) The second lowest extremal invariant
measure of the contact process, {\it Ann. Probab.} {\bf 25}, 1846--1871.

Salzano, M. \and Schonmann, R.H. (1998) The second lowest extremal invariant
measure of the contact process II, {\it Ann. Probab.}, to appear.

Schonmann, R.H. (1999a) Stability of infinite clusters in supercritical
percolation, {\it Probab. Theory Rel. Fields} {\bf 113}, 287--300.

Schonmann, R.H. (1999b) percolation in $\infty+1$ dimensions at the
uniqueness threshold, in {\it Perplexing Probability Problems: Papers in
Honor of H. Kesten}, M. Bramson and R. Durrett (eds.), Birkh\"auser, to
appear.

Trofimov, V.I (1985) Groups of automorphisms of graphs as topological
groups, {\it Math. Notes} {\bf 38}, 717--720.

Zuk, A. (1996) La propri\'et\'e (T) de Kazhdan pour les groupes agissant
sur les poly\`edres, {\it C. R. Acad. Sci. Paris S\'er. I Math.} {\bf 323},
453--458.

\endreferences

\filbreak
\begingroup
\eightpoint\sc
\parindent=0pt\baselineskip=10pt

Department of Mathematics,
Indiana University,
Bloomington, IN 47405-5701, USA
\emailwww{rdlyons@indiana.edu}
{http://php.indiana.edu/\string~rdlyons/}

Mathematics Department,
The Weizmann Institute of Science,
Rehovot 76100, Israel
\emailwww{schramm@wisdom.weizmann.ac.il}
{http://www.wisdom.weizmann.ac.il/\string~schramm/}

\endgroup

\bye